\documentclass{article}
\usepackage{graphicx} % Required for inserting images
\usepackage{amsmath,amsfonts,amssymb,amsthm}
\usepackage{tabularx}
\usepackage{cite}
\usepackage{tablefootnote}
\usepackage{bm}
\usepackage{array}
\usepackage{paracol}
\usepackage{float}
\usepackage{booktabs}     % nicer rules

\newcolumntype{C}{>{$}c<{$}}  % math mode centered column

\usepackage{geometry}
\newtheorem{theorem}{Theorem}[section]
\newtheorem{conjecture}{Conjecture}[section]
\newtheorem{definition}{Definition}[section]
\newtheorem{example}{Example}[section]

\newtheorem{remark}{Remark}[section]

\title{On the Orbits of Similarity Classes of Tetrahedra Generated by the Longest-Edge Bisection Algorithm}
\author{Jérôme Michaud and Sergey Korotov}

\date{\today}

\begin{document}
\maketitle
\begin{center}
Division of Mathematics and Physics, UKK, Mälardalen University, Västerås, Sweden \newline
emails: jerome.michaud@mdu.se,  sergey.korotov@mdu.se 
\end{center}

\begin{abstract}
 In this work, we study the dynamics of similarity classes of tetrahedra generated by the longest-edge bisection (LEB) algorithm. Building on the normalization strategy introduced by Perdomo and Plaza for triangles in \cite{Perdomo2013,Perdomo2014}, we construct a canonical representation of tetrahedra in a normalized space embedded in the product of the hyperbolic half-plane and the hyperbolic half-space model. This representation allows us to define the left and right refinement maps, $\Phi_L$ and $\Phi_R$, acting on the space of normalized tetrahedral shapes, and to study their iterative orbits as discrete dynamical systems. Using these maps, we show that the orbit of the space-filling Sommerville tetrahedron contains only 4 similarity classes, 3 of which form an attractive cycle corresponding to the orbit of the path tetrahedron. We also show that small perturbations of elements in those orbits still lead to finite orbits. In addition, we study small perturbations of the regular tetrahedron and show that their orbits are also finite. Extensive numerical exploration of orbits for the other types of tetrahedra suggests that the LEB algorithm does not produce degenerating tetrahedra.  Our framework provides a geometric and dynamical foundation for analyzing the shape evolution of tetrahedral meshes and offers a possible route toward an analytic proof of the non-degeneracy property for the tetrahedral partitions generated by the LEB refinements. This property is highly desired in e.g. the finite element methods (FEMs).
\end{abstract}

\medskip



\section{Introduction\label{sec:intro}}

When solving boundary value problems using the finite element method (FEM), the computational domain is usually divided into small,  non-overlapping, and often simplex-shaped, elements (whose total set is commonly called mesh). Ideally, we should also have a suitable refinement technique which produce smaller and smaller elements (i.e. finer meshes). In three dimensions, this means partitions of the domain of interest into tetrahedra whose shape and quality play a crucial role in the accuracy and stability of the numerical approximation \cite{shewchuk2002good}. In particular, to guarantee the optimal convergence rates and bounded interpolation constants, it is necessary that the generated families of tetrahedral partitions remain \emph{regular} or \emph{strongly regular} in the sense of Zlámal, see~\cite{Zlamal1968,Brandts2008,Brandts2011}. 

Among the various refinement techniques developed for generating adaptive and hierarchical FEM meshes, the \emph{longest-edge bisection} (LEB) algorithm has  emerged as one of the most widely used and conceptually simplest strategies. Originally introduced in the context of solving nonlinear systems of equations by Stynes~\cite{Stynes1980} and later adopted for mesh generation/refinements in various works  by Rivara, see e.g.~\cite{Rivara1984,Rivara1997}, the method recursively bisects each simplex by a plane passing through the midpoint of its longest edge and the remaining vertices. In two dimensions, the algorithm has been shown to produce a finite number of similarity classes of triangles and thus to preserve the non-degeneracy of the meshes~\cite{Rosenberg1975,Stynes1979}. However, in three dimensions, the situation is considerably more complicated. While the method idea can be straightforwardly extended to tetrahedral meshes, a general proof of regularity remains elusive for several decades~\cite{Korotov2008,Hannukainen2014}.

Recent studies on the topic have therefore concentrated on the classification of tetrahedral elements generated by repeated application of the LEB algorithm according to their \emph{similarity classes}. Each refinement step produces new tetrahedra whose shape can be characterized up to translation, rotation, reflection, and scaling. If the number of distinct similarity classes produced during refinement is finite, then regularity follows immediately. This approach has led to a number of algorithmic and computational investigations. In particular, Suárez, Trujillo-Pino, and Moreno~\cite{Suarez2021} developed an integer-arithmetic algorithm that normalizes tetrahedra by their ordered sextuple of squared edge lengths, allowing an exact comparison of similarity classes as long as squared edge lengths are integers. In general this normalization procedure is not exact, but numerically precise enough to identify similarity classes. They call their algorithm the \emph{Similarity classes longest-edge bisection} (SCLEB) algorithm. Trujillo-Pino \emph{et al.}~\cite{Trujillo2024} refined this analysis and proved that nearly equilateral tetrahedra generate at most thirty-seven distinct similarity classes under successive bisections. Other related works have examined special cases, such as the regular or trirectangular tetrahedra, within the framework of the eight-tetrahedron partition~\cite{Padron2023}. These approaches have significantly advanced the enumeration of similarity classes but remain largely combinatorial in nature. They rely on edge-length representations and discrete comparison criteria, which do not easily lend themselves to analytical study of the underlying geometric transformations induced by the LEB process.

In contrast, the situation in two dimensions has benefited from a \emph{geometric normalization} of triangular elements. Perdomo and Plaza~\cite{Perdomo2013,Perdomo2014} introduced a space of normalized triangles endowed with a hyperbolic metric in the Poincaré half-plane. In this framework, each triangle corresponds to a point in a canonical region of the hyperbolic plane, and the refinement operation acts as an isometry of that space. This geometric viewpoint allows for both a continuous description of the discrete refinement process and an analytic proof of non-degeneracy for the resulting families of triangulations.

The aim of the present work is to extend this geometric perspective to the three-dimensional setting. We introduce a normalization procedure that maps each tetrahedral element to a canonical representative in a space of \emph{normalized tetrahedra} obtained by selecting a base face to be mapped onto Perdomo and Plaza's space of triangle \cite{Perdomo2013,Perdomo2014} and the fourth vertex to the upper half-space using shape preserving transformations. The resulting representation of a tetrahedron can be seen as a point in the product space of the Poincaré half-plane $\mathbb H^2$ and the Poincaré half-space $\mathbb H^3$ models. Using the induced product metric, we obtain a canonical metric space of tetrahedra. Within this space, the LEB refinement step is represented by two maps, $\Phi_L$ and $\Phi_R$, corresponding to the left and right sub-tetrahedra generated by bisecting the longest edge. The iterative behavior of these maps defines orbits of normalized tetrahedra whose closure characterizes the long-term dynamics of the refinement process. By studying these orbits, we seek to identify closed families of tetrahedra, describe their orbits, and provide numerical evidence of boundedness and non-degeneracy of the resulting partitions.

This normalization framework thus connects the algorithmic classification of similarity classes of Suárez, Trujillo-Pino, and Moreno \cite{Suarez2021} with a geometric and dynamical-systems interpretation of the refinement process. It offers a unified setting in which both analytical and numerical investigations of tetrahedral regularity can be pursued.

\section{Normalizing and visualizing the similarity classes of tetrahedra\label{sec:norm}}

In this section, we introduce a normalization procedure inspired by the work of Perdomo and Plaza \cite{Perdomo2013,Perdomo2014} for triangles. The idea is to select one of six faces of the tetrahedron as its canonical base and send that face to the canonical space of triangles by translation, rotations, and scaling operations, hence preserving the shape of the tetrahedron. When doing such a transformation, the fourth vertex can have either a positive or negative third component. Therefore, similarly to \cite{Perdomo2013,Perdomo2014}, we also use a reflection (if needed) to ensure that the fourth vertex is finally mapped to the upper half-space. Using such a normalization procedure, we can study the longest-edge bisection (LEB) refinement process both analytically and numerically as a discrete dynamical system. Furthermore, such type of normalization provides a natural way to visualize the orbit of a given tetrahedron under the LEB refinement process. Compared to other normalization approaches, such as that of \cite{Suarez2021,Trujillo2024} based on the square length of edges, our normalization procedure uniquely maps similarity classes to a compact metric space and provides a natural and intuitive visualization of the underlying dynamics.

%In this section we describe the geometric framework used to study the  longest-edge bisection (LEB) of tetrahedra through normalization.  The goal is to represent every tetrahedral element by a unique canonical  representative belonging to a compact space of normalized shapes.  This allows us to analyze the refinement process as the iteration of  well-defined maps acting on this space, thus transforming the discrete  refinement algorithm into a dynamical system.

\subsection{The canonical space of tetrahedra}
%Let $A,B,C,D \in \mathbb{R}^3$ be the vertices of a non-degenerate tetrahedron $ABCD$. %By non-degenerate, we mean that the volume of $ABCD$ is strictly positive. 

\begin{definition}
Let $\mathcal{S}$ be the space of all non-degenerate tetrahedra.
    Two tetrahedra $T_1\in \mathcal{S}$ and $T_2\in \mathcal{S}$ are called \emph{shape similar} if and only if there exists an isometry $R\in O(3)$, a translation vector $b$, and a positive scalar $\lambda >0$, such that $T_2 =b + \lambda R T_1$, and we write $T_1 \sim T_2$. The space of similarity classes is then defined as $\mathcal{S}/\sim$.
\end{definition}

The normalization procedure is, in general, about selecting a suitable (canonical) representative of each class in $\mathcal{S}/\sim$. In order to describe that in detail, we need to introduce some extra definitions and denotation.

\medskip

First, we recall the concept of the so-called canonical space of triangles, see \cite{Perdomo2013}.

\begin{definition}
    The canonical space of triangles consists of triangles having $(0,0)$ and $(0,1)$ as two (fixed) vertices in the complex plane, and for which the third vertex lies in the set  
    \[
\mathcal B := \{z \in \mathbb C\ |\  \Im(z)>0,\ 0<\Re(z)\leq\frac12,\ |z-1|\leq 1\}\,,
\]
where $\Im(z)$ is the imaginary part of $z$ and $\Re(z)$ is the real part of $z$.
\end{definition}

\begin{figure}
    \centering
    \includegraphics[width=0.5\linewidth]{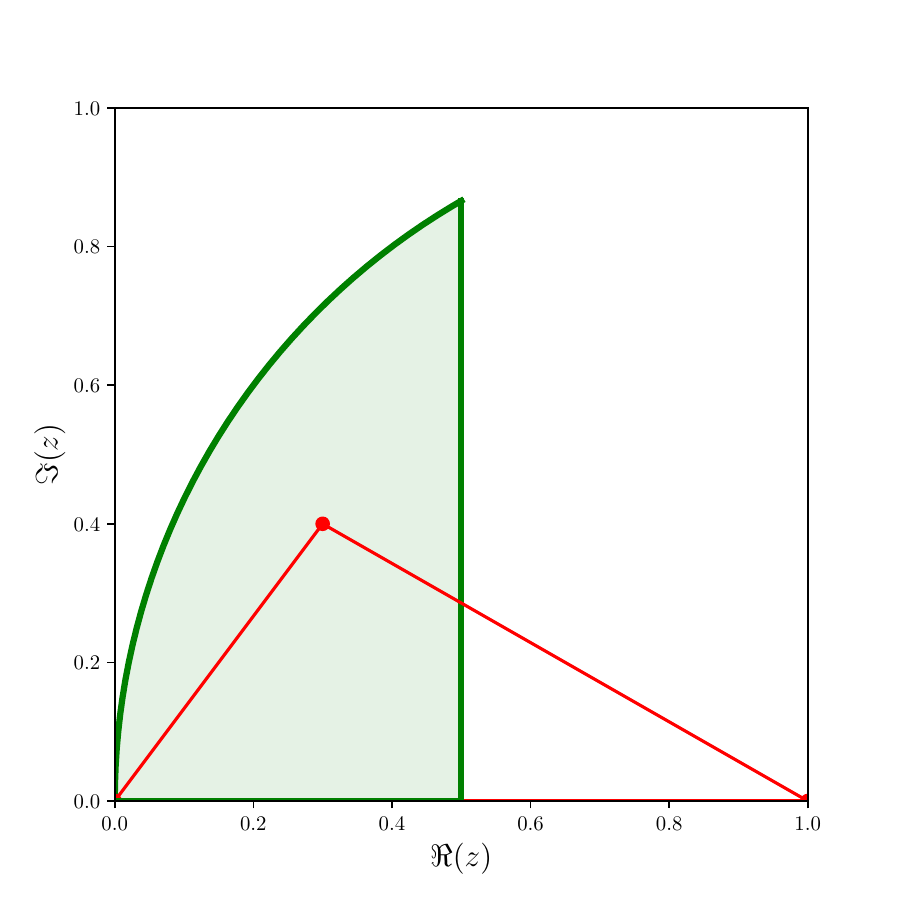}
    \caption{An illustration for the canonical space of triangles. The set $\mathcal B$, of admissible positions for the third vertex, is colored in green.}
    \label{fig:space_of_triangles}
\end{figure}

\begin{remark}
  From the above definition, it follows that the edge connecting $0$ and $1$ (of the length 1) is the longest (or one of the longest ones), and the edge connecting $0$ and $z$, of the length  $|z|$, is the shortest (or one of the shortest ones) for all triangles from the canonical space of triangles, see Figure~\ref{fig:space_of_triangles}.  
\end{remark}

\medskip

Now, each canonical tetrahedron will be defined as that having one of its faces in the canonical space of triangles, which means that three of its vertices are the points  $0$, $1$, and $z$. The set of possible positions for the fourth vertex $\mathcal V(z)$, provided $z$ is fixed, is described in the next definition.

\begin{definition}
The set of possible positions of the fourth-vertex $\mathcal V(z)$ is defined as
    \[
\mathcal V(z) := \{(w,t)\in \mathbb{C}\times \mathbb R_+\ |\ |z|^2\leq |w|^2 + t^2 \leq 1,\ |z|^2\leq|w-1|^2 + t^2 \leq 1 ,\ |w-z|^2 + t^2 \leq 1\}\,.
\]
\end{definition}
If $z= z_1+z_2i$ and $w=w_1 + w_2i$, then each canonical tetrahedron is uniquely associated with a tetrahedron with vertices  $(0,0,0)$, $(1,0,0)$, $(z_1,z_2,0)$, and $(w_1,w_2,t)$.

\begin{definition}
The canonical space of tetrahedra $\mathcal{T}$ is defined as
\[
\mathcal{T}:=\{(z,(w,t))\in \mathcal{B}\times \mathcal{V}(z)\}\,.
\]
Elements of $\mathcal T$ are denoted by $\tau$, possibly with indices.
\end{definition}

\begin{figure}
    \centering
    \includegraphics[width=0.45\linewidth]{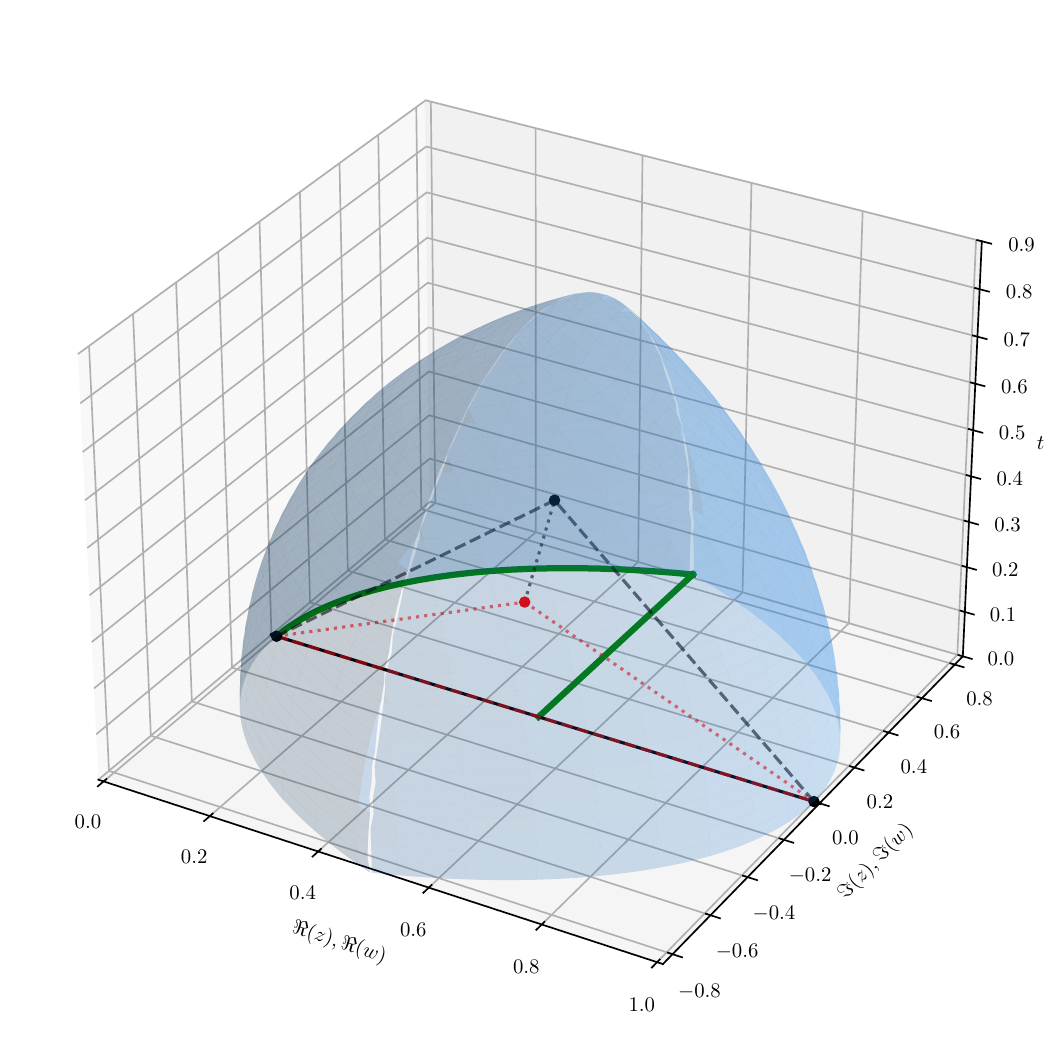} \hfill \includegraphics[width=0.45\linewidth]{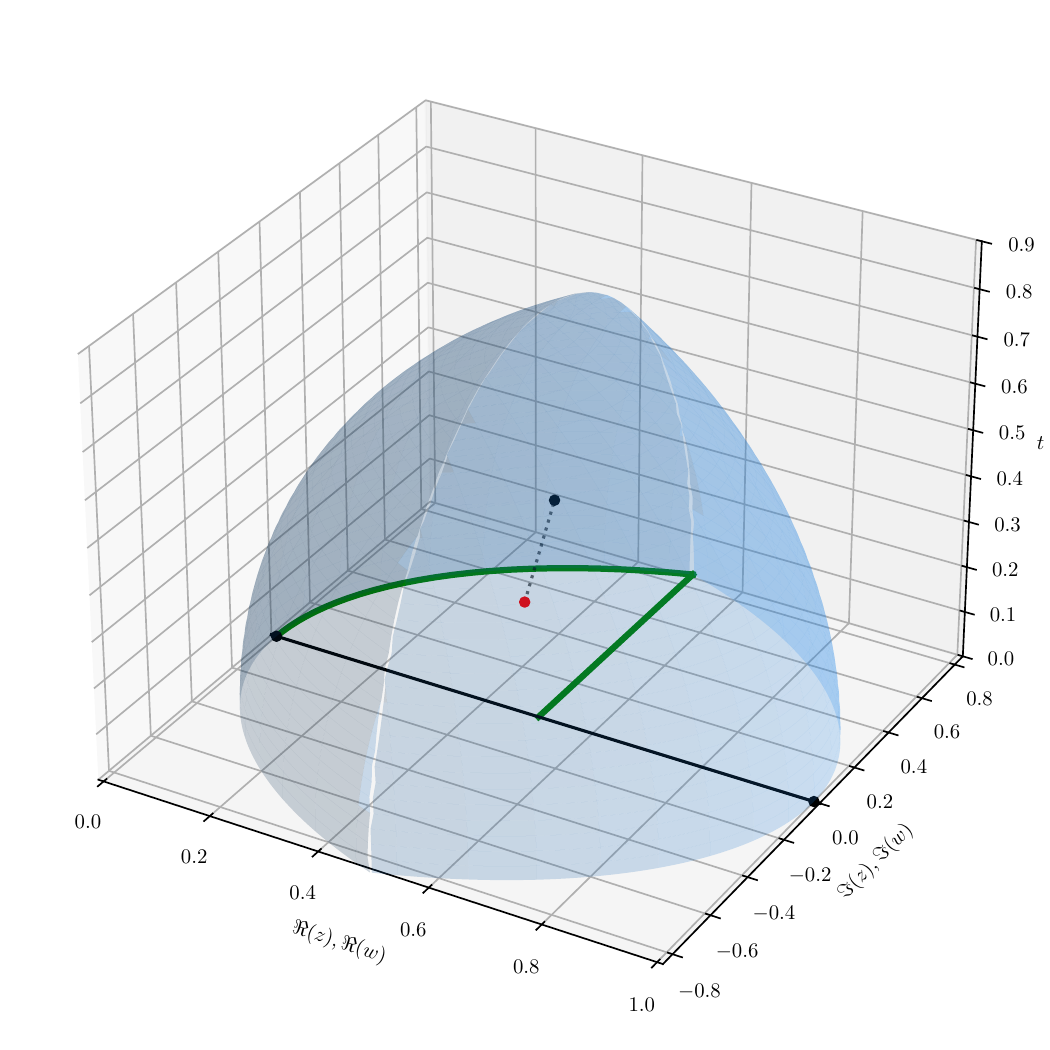} 
    \caption{Left: an example of a canonical tetrahedron. The red point is the normalized representation of the canonical base and the black point is the remaining vertex. All the edges of the tetrahedron are also displayed. Right: same figure, but only the edge carrying information is presented. }
    \label{fig:viz}
\end{figure}

An example of visualization for canonical tetrahedra generalizing the case of canonical triangles from \cite{Perdomo2013} is illustrated in Figure \ref{fig:viz}. In addition to the space of triangles marked by the green color, we added the intersection of two spheres of radius $1$ centered at $(0,0,0)$ and at $(1,0,0)$ in the upper half-space. Such a visual representation of all tetrahedra generated by LEB will be used in the rest of this paper.

\medskip

In the spirit of \cite{Perdomo2013,Perdomo2014}, we can endow $\mathcal T$ with a metric by using the fact that $\mathcal T \subset \mathbb{H}^2 \times \mathbb{H}^3$ and defining the product metric.
\begin{definition}\label{def:metric}
    Let $\mathcal{T} \subset \mathbb{H}^2 \times \mathbb{H}^3$, 
with points $(z,(w,t)), \; z \in \mathcal B \subset \mathbb{H}^2, \; (w,t) \in \mathcal V(z) \subset \mathbb{H}^3$.
The product metric is:
\[
ds^2 = \frac{|dz|^2}{(\Im (z))^2} + \frac{|dw|^2 + dt^2}{t^2}\,.
\]
\end{definition}

\begin{definition}\label{def:distance}
The distance between two points 
$\tau_1 =(z_1,(w_1,t_1)), \tau_2=(z_2,(w_2,t_2)) \in \mathcal{T}$
 is given by
\[
d_{\mathcal{T}}\big(\tau_1,\tau_2\big)
= \sqrt{\, d_{\mathbb{H}^2}(z_1,z_2)^2 + d_{\mathbb{H}^3}((w_1,t_1),(w_2,t_2))^2 \,}.
\]
Here the hyperbolic distances are:
\[
\cosh d_{\mathbb{H}^2}(z_1,z_2)
= 1 + \frac{|z_1 - z_2|^2}{2\,\Im (z_1)\,\Im (z_2)},
\]
\[
\cosh d_{\mathbb{H}^3}((w_1,t_1),(w_2,t_2))
= \frac{|w_1 - w_2|^2 + t_1^2 + t_2^2}{2\,t_1\,t_2}.
\]
\end{definition}

%Now, we want to design a normalization procedure that ensures that there is a bijection between $\mathcal{S}/\sim$ and $\mathcal{T}$.  

\subsection{Normalization procedure}
Let $T\in \mathcal{S}$ be a non-degenerate tetrahedron. Our normalization procedure maps a well-chosen face with vertices denoted by $A,B,C$ to $\mathcal B$ and the vertex $D$ to $\mathcal V(z)$, hence, mapping $T$ to an element $\tau$ of $\mathcal{T}$. 

\medskip
\begin{definition}
  The normalization procedure consists of the following steps:

\begin{enumerate}
\item \textbf{Selection of the base face and labeling of vertices.}
  Let us assume that the edge $AB$ is the longest edge of  $T$, and let $AC$ be the shortest among four edges adjacent to $AB$. If the ambiguity is still not resolved, the face with the smallest angle $\angle{BAC}$ is chosen. If the ambiguity is still not resolved, the choice is arbitrary and will lead to the same normalization. The face $ABC$ is then selected as the base.

\item \textbf{Rigid motion and scaling.}  
  Apply a translation, rotation and scaling so that $A$ is mapped to $(0,0,0)$ and $B$ to $(1,0,0)$, and the base face $ABC$ lies in the plane $\{x_3=0\}$ with $C$ in the upper half-plane ($x_2>0$). In practice, let us define $\bm{a}, \bm{b}, \bm{c}$, and $\bm{d}$ the vectors corresponding to the vertices $A,B,C$ and $D$. 
  \begin{enumerate}
      \item   Translate $A$ to the origin. The coordinates of the other vertices are then defined by $\bm{b}' = \bm{b}-\bm{a}$, $\bm{c}' = \bm{c}-\bm{a}$, $\bm{d}' = \bm{d}-\bm{a}$.
      \item Construct the orthonormal frame $\bm{e}_1 = \frac{\bm{b}'}{\|\bm{b}'\|}$, $\bm{e}_2 = \frac{\bm{v}_2}{\|\bm{v}_2\|}$, where $\bm{v}_2 = \bm{c}' - (\bm{c}' \cdot \bm{e}_1) \bm{e}_1$ and $\bm{e}_3 = \bm{e}_1 \times \bm{e}_2$.
      \item Rotate the tetrahedron using $\bm{b}'' = R\bm{b}'$, $\bm{c}'' = R\bm{c}'$, $\bm{d}'' = R\bm{d}'$, where $R = \begin{bmatrix}
\bm{e}_1^\top \\
\bm{e}_2^\top \\
\bm{e}_3^\top
\end{bmatrix}$.
\item Scale the tetrahedron by defining $\bm{b}''' = \frac{\bm{b}''}{\|\bm{b}''\|}$, $\bm{c}''' = \frac{\bm{c}''}{\|\bm{b}''\|}$, $\bm{d}''' = \frac{\bm{d}''}{\|\bm{b}''\|}$.
  \end{enumerate}
  The coordinates of $C$ are thus of the form  $C=(c'''_1,c'''_2,0)$. This corresponds to a point in the space of triangles $z = c'''_1 +ic'''_2 \in \mathcal{B}$.

\item \textbf{Placement of the fourth vertex.}  
  After step 2. The coordinates of $D$ are $(d'''_1,d'''_2,d'''_3)$. To make sure that the point $D$ ends up in the upper half-space, reflect if needed the third component to the upper half-space. In practice, $\bm{d}''''=(d'''_1,d'''_2,|d'''_3|)$. Since the tetrahedron $\tau$ is non-degenerate, we get $d''''_3>0$. This corresponds to a point $(w,t)\in \mathcal{V}(z)$, with $w:=d'''_1 +d'''_2i$ and $t=d''''_3$.
\end{enumerate}

The resulting configuration $(z,(w,t))\in \mathcal{T} $ defines the normalized form of the tetrahedron $T$.  
\end{definition}
 
\begin{theorem}
    The normalization map
\[
\begin{aligned}
    N : \mathcal{S} &\longrightarrow \mathcal{T}\\
    T &\longmapsto \tau = (z,(w,t))
\end{aligned}
\]
is surjective and constant on similarity classes. Moreover, for two tetrahedra $T_1$ and $T_2$ that are not shape similar, $T_1 \not\sim T_2$, then $N(T_1)\neq N(T_2)$.
\end{theorem}
\begin{proof}
    By definition of the normalization procedure, any element $T\in \mathcal S$ is sent to $\mathcal T$, proving the surjectivity of $N$. When the longest canonical face is uniquely defined by step 1 of the normalization procedure, it follows that two shape similar tetrahedra will be mapped to the same element of $\mathcal T$. 
    
    The only case that need some care is when step 1 is still ambiguous because either two or more faces are similar or when the canonical face is isosceles with the two equal edges being the longest. In both cases, either step 1 and 2 in the normalization procedure lead to the same representation, or they lead to the vertex $D$ being map either to the upper half-space, or to the lower half-space. Step 3 of the normalization procedure ensures that the choice of the canonical face is arbitrary and that the normalized form $N(T)$ is the same. 

    Finally, since two elements $\tau_1\neq\tau_2$ of $\mathcal T$ are not shape similar, it follows that two tetraheda $T_1\not\sim T_2$ will be mapped to different elements of $\mathcal T$, since shape similarity is preserved by the normalization map $N$.
\end{proof}

\begin{example}\label{ex:1}
    Let us consider the space-filling Sommerville tetrahedron \cite{Sommerville1923} with vertices $(-1,0,0)$, $(1,0,0)$, $(0,-1,1)$, and $(0,1,1)$.  Applying the normalization procedure to this tetrahedron leads to the normalized representation $\tau^S:=\left(\frac12+ \frac{\sqrt{2}}{2}i,\left(\frac12+0i,\frac{\sqrt{2}}{2}\right)\right)$.
\end{example}

\begin{example}\label{ex:2}
    Let us consider the path tetrahedron with vertices $(0,0,0)$, $(1,0,0)$, $(1,1,0)$, and $(1,1,1)$. Applying the normalization procedure to this tetrahedron leads to the normalized representation $\tau^P:=\left(\frac13+ \frac{\sqrt{2}}{3}i,\left(\frac23+\frac{\sqrt{2}}{6}i,\frac{\sqrt{6}}{6}\right)\right)$.
    \end{example}

Compared to the normalization procedure based on the square of edge lengths \cite{Suarez2021,Trujillo2024}, our normalization procedure yields an exact 5-dimensional representation of tetrahedra. In the rest of this paper, the normalization procedure will be carried out symbolically using the SymPy package \cite{SymPy} to avoid numerical inaccuracies and ensure the correct identification of similarity classes mirroring the integer-based strategy used in \cite{Suarez2021,Trujillo2024}.

\begin{remark}\label{rem:norm}
    It is possible to modify our methodology to resemble that of \cite{Suarez2021} by changing the base face used in the normalization procedure. In \cite{Suarez2021}, the canonical face is defined by the longest edge, the longest adjacent edge to the longest. Ties are resolved by choosing the face with the longest edge closing the face. The difference in choices of canonical faces used during LEB refinements can lead to orbits of different lengths. However, both strategies are, obviously, still in the general framework of the LEB strategy.
\end{remark}

\subsection{LEB in the space of normalized tetrahedra}
Let $T\in \mathcal{S}$ be a non-degenerate tetrahedron, $N(T) =(z,(w,t))$ its normalization and $\hat{T}=ABCD\in \mathcal{S}$ the corresponding canonical tetrahedron with vertices $A=(0,0,0)$, $B=(1,0,0)$, $C=(z_1,z_2,0)$ and $D=(w_1,w_2,t)$. Under the LEB refinement, the tetrahedron $\hat{T}$ is bisected by the plane passing through the midpoint $M$ of its longest edge $AB$ and the opposite vertices $C$ and $D$. Since the longest edge of $\hat{T}$ is $AB$, we define the vertex $M=(0.5,0,0)$ at the mid-point of $AB$. The left and right sub-tetrahedra $T_L$ and $T_R$ of $\hat{T}$ are defined then as $T_L=AMCD$ and $T_R = MBCD$. 

\begin{definition}
    Let $T\in\mathcal{S}$ be a non-degenerate tetrahedron and $\tau= N(T)$ its normalized representation. The left and right maps associated to the LEB process are defined as the following:
    \[
\begin{aligned}
    \Phi_{L,R} : \mathcal{T} &\longrightarrow \mathcal{T}\\
    \tau&\longmapsto \tau_{L,R} = N(T_{L,R} )
\end{aligned}\,,
\]
respectively.
\end{definition}

\begin{remark}
    The maps $\Phi_L$ and $\Phi_R$ can be explicitly defined as functions of $\tau$ in a similar manner to what is done in \cite{Perdomo2013,Perdomo2014}, but we are not doing it here. To do it, we would need to write down a formula for each of the 24 possible vertex orderings defined in step 1 of the normalization process for both the left and the right map, leading to 48 different expressions. In addition, we should specify conditions under which each transformation applies. 
\end{remark}

\begin{definition}
    A region $\Omega \subset \mathcal{T}$ is called a \emph{closed region for LEB} if for all $\tau \in \Omega$, we have  $\Phi_L(\tau) \in \Omega$ and $\Phi_R(\tau)\in \Omega$.
\end{definition}

\begin{remark}
    By definition of $\Phi_{L,R}$, $\mathcal T$ is a closed region for LEB.
\end{remark}

\begin{definition}
    Let $\tau\in \mathcal{T}$. We define $\Gamma_{\tau}^{(0)} = \{\tau\}$,  and $\Gamma_{\tau}^{(n+1)} = \Phi_L(\Gamma_\tau^{(n)}) \cup \Phi_R(\Gamma_\tau^{(n)})$ for $n\geq 0$. The \emph{orbit} of $\tau$ is defined then as the set $\Gamma_\tau:= \bigcup_{n\geq0}\Gamma_\tau^{(n)}$.
\end{definition}

For any $\tau\in \mathcal T$, its orbit $\Gamma_{\tau}$ is closed by definition. In the rest of this paper, we will study the orbits of some specific tetrahedra.

\section{On orbits produced by the LEB refinement algorithms\label{sec:orbits}}
In this section, we study the orbits of various (normalized) tetrahedra under the LEB refinement algorithm. We will show that some tetrahedra have very short orbits. As pointed out already in \cite{Hannukainen2014}, the orbit of the Sommerville tetrahedron contains only four elements and the orbit of the path tetrahedron contains only three elements. This has also been noted by Trujillo-Pino \emph{et al.} \cite{Trujillo2024}, where the orbit of the Sommerville tetrahedron has been studied using the squared edge length normalization procedure. Here, we prove that the orbit of the path tetrahedron is actually included in that of the Sommerville tetrahedron. In addition, we show that small perturbations of elements from these orbits still lead to finite but longer orbits. The we use our algorithm to study nearly regular tetrahedra and compare it to previous results \cite{Suarez2021}. Finally, we  study four cases that do not seem to have finite orbits.

\subsection{The case of the Sommerville tetrahedron}
The space-filling Sommerville tetrahedron \cite{Sommerville1923} plays a special role in the LEB refinement algorithm as well as in numerical analysis in general\cite{hovsek2017role}. The following result holds.

\begin{figure}[ht]
    \centering
    \includegraphics[width=0.5\linewidth]{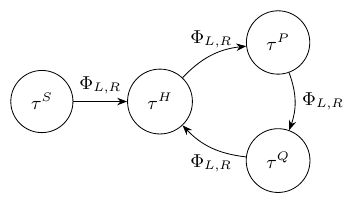}
    \caption{The orbit of the Sommerville tetrahedron.}
    \label{fig:graph_som}
\end{figure}

\begin{theorem}\label{prop:Sommerville}
    The orbit of the Sommerville tetrahedron contains four elements. Its orbit is attracted into a cycle of length three including the path tethrahedron. Hence, the orbit of the path tetrahedron contains exactly three elements. 
\end{theorem}
\begin{proof}
    As mentioned in Example \ref{ex:1}, the normalized representation of the Sommerville tetrahedron is $\tau^S:=\left(\frac12+ \frac{\sqrt{2}i}{2},\left(\frac12+0i,\frac{\sqrt{2}}{2}\right)\right)$. Applying the left and right maps to $\tau^S$ and iterating, we find that
    \begin{itemize}
        \item $\Phi_L(\tau^S)=\Phi_R(\tau^S) = \left(\frac12+\frac{i}2,\left(\frac12+\frac{i}2,\frac12\right)\right)=:\tau^H$,
        \item $\Phi_L(\tau_1)=\Phi_R(\tau_1) = \left(\frac13+ \frac{\sqrt{2}i}{3},\left(\frac23+\frac{\sqrt{2}i}{6},\frac{\sqrt{6}}{6}\right)\right)=:\tau^P$,
        \item $\Phi_L(\tau^P)=\Phi_R(\tau^P) = \left(\frac12+ \frac{\sqrt{2}i}{4},\left(\frac12+0i,\frac12\right)\right)=:\tau^Q$,
        \item $\Phi_L(\tau_2)=\Phi_R(\tau_2)=\tau^H$.

    \end{itemize}
    The proof is illustrated in Figure~\ref{fig:graph_som}.
\end{proof}
\begin{remark}\label{rem:names}
    In order to simplify further discussion, we denote each element in the orbit of the Sommerville tetrahedron $\tau^S$ as marked in  Figure \ref{fig:graph_som}. We use the denotation  $\tau^H$ for the first element because all its non-trivial components are equal to $\frac12$ and $H$ stands for \emph{halves}. $\tau^P$ is the path tetrahedron and $\tau^Q$ is its image under $\Phi_{L,R}$.
\end{remark}

\begin{figure}
    \centering
    \includegraphics[width=0.85\linewidth]{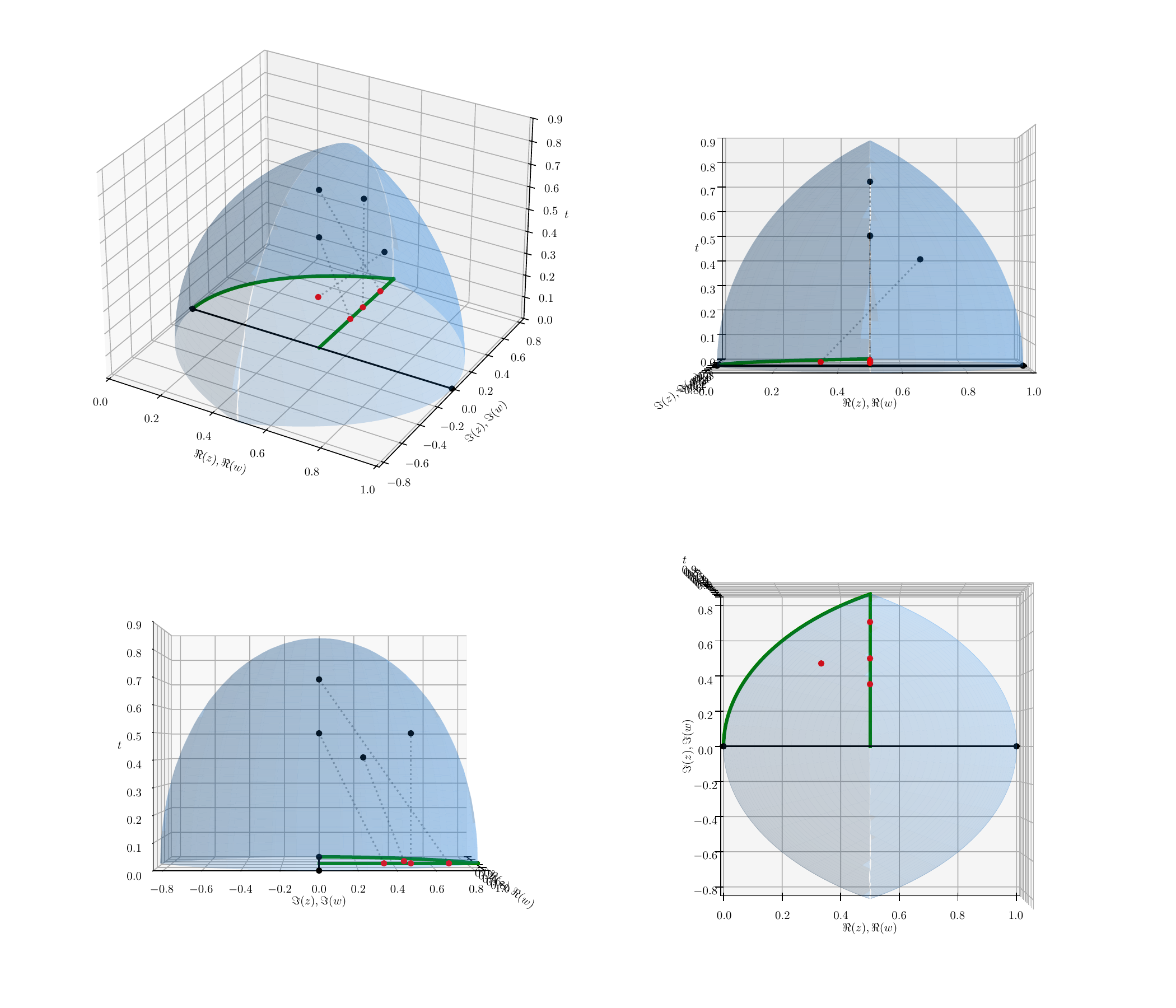}
    \caption{Visualization of the orbit of the Sommerville tetrahedron. The first three panels show different views of the same results. The bottom right panel shows only the positions of the base vertices.}
    \label{fig:Sommerville-orbit}
\end{figure}

Theorem~\ref{prop:Sommerville} confirms the result of \cite{Hannukainen2014} and \cite{Trujillo2024} for the number of shapes in the orbits of these two specific tetrahedra. Visualization of  the orbit of the Sommerville tetrahedron is presented in Figure \ref{fig:Sommerville-orbit}.

\subsection{Orbits of small perturbations of the Sommerville and path tetrahedra}
\begin{figure}[H]
    \centering
    \small
\parbox[t]{\textwidth}{
\begin{tabularx}{0.48\textwidth}{@{}lX@{}}
\hline
ID & Normalization \\
\hline
$\tau_{1}$ &  \begin{minipage}{\linewidth}$\left(\frac{190}{401} + \frac{10 \sqrt{842}}{401} i, \left(\frac{190}{401} + \frac{200 \sqrt{842}}{168821} i, \frac{210 \sqrt{337642}}{168821}\right)\right)$\end{minipage} \\
\hline
$\tau_{2}$ &  \begin{minipage}{\linewidth}$\left(\frac{1}{2} + \frac{21}{40} i, \left(\frac{1}{2} + \frac{1}{2} i, \frac{1}{2}\right)\right)$\end{minipage} \\
$\tau_{3}$ &  \begin{minipage}{\linewidth}$\left(\frac{1}{2} + \frac{10}{21} i, \left(\frac{10}{21} + \frac{10}{21} i, \frac{10}{21}\right)\right)$\end{minipage} \\
$\tau_{4}$ &  \begin{minipage}{\linewidth}$\left(\frac{1}{2} + \frac{\sqrt{401}}{40} i, \left(\frac{1}{2} + \frac{211 \sqrt{401}}{8020} i, \frac{21 \sqrt{401}}{802}\right)\right)$\end{minipage} \\
$\tau_{5}$ &  \begin{minipage}{\linewidth}$\left(\frac{1}{2} + \frac{10}{21} i, \left(\frac{11}{21} + \frac{10}{21} i, \frac{10}{21}\right)\right)$\end{minipage} \\
$\tau_{6}$ &  \begin{minipage}{\linewidth}$\left(\frac{1}{2} + \frac{10 \sqrt{401}}{401} i, \left(\frac{211}{401} + \frac{10 \sqrt{401}}{401} i, \frac{210}{401}\right)\right)$\end{minipage} \\
$\tau_{7}$ &  \begin{minipage}{\linewidth}$\left(\frac{190}{401} + \frac{210}{401} i, \left(\frac{190}{401} + \frac{210}{401} i, \frac{10 \sqrt{401}}{401}\right)\right)$\end{minipage} \\
$\tau_{8}$ &  \begin{minipage}{\linewidth}$\left(\frac{10}{21} + \frac{10}{21} i, \left(\frac{10}{21} + \frac{10}{21} i, \frac{10}{21}\right)\right)$\end{minipage} \\
$\tau_{9}$ &  \begin{minipage}{\linewidth}$\left(\frac{1}{2} + \frac{\sqrt{401}}{40} i, \left(\frac{1}{2} + \frac{19 \sqrt{401}}{802} i, \frac{21 \sqrt{401}}{802}\right)\right)$\end{minipage} \\
$\tau_{10}$ &  \begin{minipage}{\linewidth}$\left(\frac{1}{2} + \frac{10 \sqrt{401}}{401} i, \left(\frac{190}{401} + \frac{10 \sqrt{401}}{401} i, \frac{210}{401}\right)\right)$\end{minipage} \\
$\tau_{11}$ &  \begin{minipage}{\linewidth}$\left(\frac{1}{2} + \frac{21}{40} i, \left(\frac{1}{2} + \frac{11}{20} i, \frac{1}{2}\right)\right)$\end{minipage} \\
\hline
$\tau_{12}$ &  \begin{minipage}{\linewidth}$\left(\frac{1}{3} + \frac{\sqrt{2}}{3} i, \left(\frac{13}{20} + \frac{7 \sqrt{2}}{40} i, \frac{7 \sqrt{6}}{40}\right)\right)$\end{minipage} \\
$\tau_{13}$ &  \begin{minipage}{\linewidth}$\left(\frac{1}{3} + \frac{\sqrt{2}}{3} i, \left(\frac{41}{60} + \frac{19 \sqrt{2}}{120} i, \frac{7 \sqrt{6}}{40}\right)\right)$\end{minipage} \\
\hline
\end{tabularx}
\hfill
\begin{tabularx}{0.48\textwidth}{@{}lX@{}}
\hline
ID & Normalization \\
\hline
$\tau_{14}$ &  \begin{minipage}{\linewidth}$\left(\frac{19}{60} + \frac{\sqrt{842}}{60} i, \left(\frac{13}{20} + \frac{77 \sqrt{842}}{8420} i, \frac{7 \sqrt{2526}}{842}\right)\right)$\end{minipage} \\
$\tau_{15}$ &  \begin{minipage}{\linewidth}$\left(\frac{100}{321} + \frac{10 \sqrt{221}}{321} i, \left(\frac{431}{642} + \frac{1055 \sqrt{221}}{70941} i, \frac{35 \sqrt{70941}}{23647}\right)\right)$\end{minipage} \\
$\tau_{16}$ &  \begin{minipage}{\linewidth}$\left(\frac{100}{321} + \frac{10 \sqrt{221}}{321} i, \left(\frac{137}{214} + \frac{385 \sqrt{221}}{23647} i, \frac{35 \sqrt{70941}}{23647}\right)\right)$\end{minipage} \\
$\tau_{17}$ &  \begin{minipage}{\linewidth}$\left(\frac{211}{642} + \frac{5 \sqrt{842}}{321} i, \left(\frac{137}{214} + \frac{350 \sqrt{842}}{45047} i, \frac{35 \sqrt{270282}}{45047}\right)\right)$\end{minipage} \\
\hline
$\tau_{18}$ &  \begin{minipage}{\linewidth}$\left(\frac{410}{841} + \frac{10 \sqrt{842}}{841} i, \left(\frac{441}{841}  - \frac{210 \sqrt{842}}{354061} i, \frac{210 \sqrt{842}}{12209}\right)\right)$\end{minipage} \\
$\tau_{19}$ &  \begin{minipage}{\linewidth}$\left(\frac{130}{267} + \frac{70 \sqrt{2}}{267} i, \left(\frac{401}{801}  - \frac{10 \sqrt{2}}{801} i, \frac{10 \sqrt{178}}{267}\right)\right)$\end{minipage} \\
$\tau_{20}$ &  \begin{minipage}{\linewidth}$\left(\frac{410}{841} + \frac{10 \sqrt{842}}{841} i, \left(\frac{400}{841} + \frac{210 \sqrt{842}}{354061} i, \frac{210 \sqrt{842}}{12209}\right)\right)$\end{minipage} \\
$\tau_{21}$ &  \begin{minipage}{\linewidth}$\left(\frac{1}{2} + \frac{\sqrt{2}}{4} i, \left(\frac{19}{40} + 0 i, \frac{21}{40}\right)\right)$\end{minipage} \\
$\tau_{22}$ &  \begin{minipage}{\linewidth}$\left(\frac{1}{2} + \frac{\sqrt{2}}{4} i, \left(\frac{21}{40} + 0 i, \frac{21}{40}\right)\right)$\end{minipage} \\
$\tau_{23}$ &  \begin{minipage}{\linewidth}$\left(\frac{130}{267} + \frac{70 \sqrt{2}}{267} i, \left(\frac{400}{801} + \frac{10 \sqrt{2}}{801} i, \frac{10 \sqrt{178}}{267}\right)\right)$\end{minipage} \\
$\tau_{24}$ &  \begin{minipage}{\linewidth}$\left(\frac{1}{2} + \frac{5 \sqrt{221}}{221} i, \left(\frac{231}{442} + 0 i, \frac{105}{221}\right)\right)$\end{minipage} \\
$\tau_{25}$ &  \begin{minipage}{\linewidth}$\left(\frac{1}{2} + \frac{5 \sqrt{221}}{221} i, \left(\frac{211}{442} + 0 i, \frac{105}{221}\right)\right)$\end{minipage} \\
\hline
\end{tabularx}}

    \includegraphics[width=0.95\linewidth]{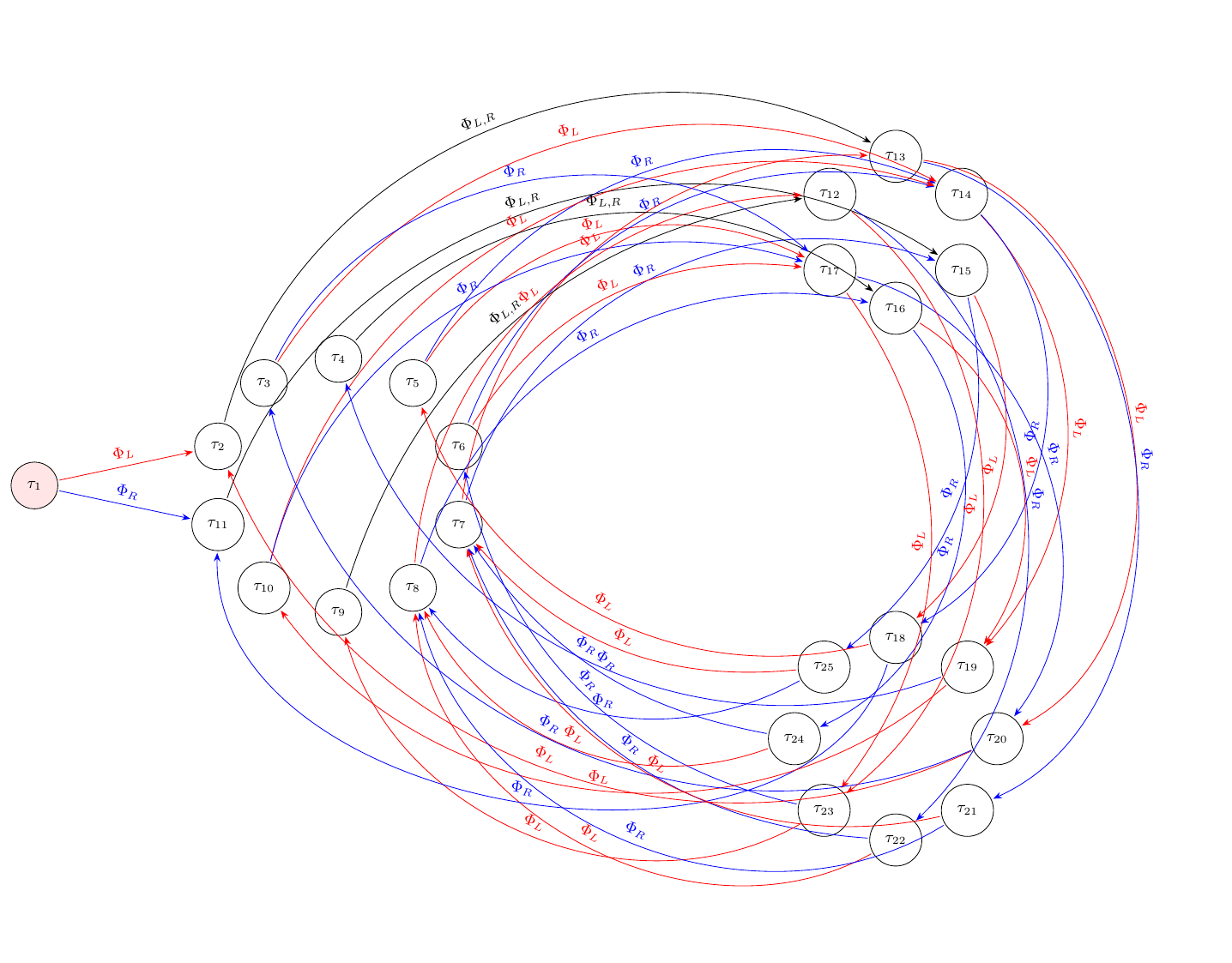}
    \caption{The results of Example~\ref{ex:Som-mod} illustrated. Graph shows how the normalized tetrahedra listed in the tables at the top and marked by nodes in the diagram map onto each other under the action of $\Phi_L$ (red) and $\Phi_R$ (blue). All the tetrahedra involved are grouped according to how they cluster. The red node $\tau_1$  corresponds to the modification of the Sommerville tetrahedron, and it has no incoming edges.}
    \label{fig:graph_Sommerville_mod}
\end{figure}

\begin{figure}[ht]
    \centering
    \includegraphics[width=0.85\linewidth]{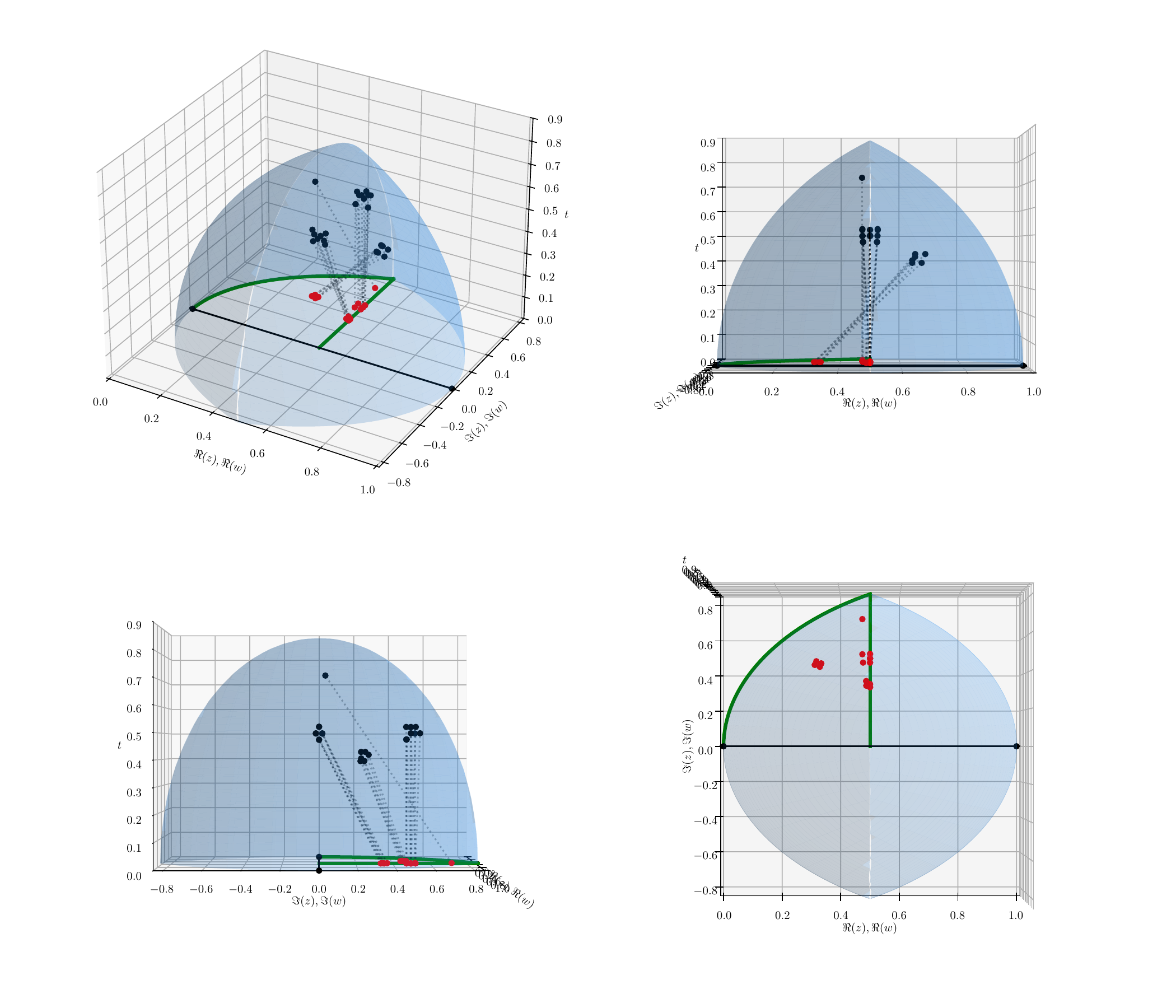}
    \caption{ Orbit of the perturbed Sommerville tetrahedra with vertices $(-1,0,0)$, $(1,0,0)$, $(0,-1,1)$, and $(0,1,\frac{11}{10})$. The bottom right panel shows only the positions of the base vertices.}
    \label{fig:Orbit_Sommerville_mod}
\end{figure}

In \cite{Perdomo2014}, the idea for constructing the proof of regularity of the bisection process essentially consists of introducing the hyperbolic metric and showing that some regions of the canonical space of triangles are closed. In order to test whether a certain neighborhood of the orbit of the Sommerville tetrahedron is closed, we make small perturbations of the Sommerville tetrahedron and of some elements from its orbit and monitor the resulting orbits.
\begin{example}\label{ex:Som-mod} We apply LEB refinements to the tetrahedron with vertices $(-1,0,0)$, $(1,0,0)$, $(0,-1,1)$, and $(0,1,\frac{11}{10})$, which normalizes to $\tilde\tau^S=\left(\frac{190}{401} + \frac{10 \sqrt{842}}{401} i, \left(\frac{190}{401} + \frac{200 \sqrt{842}}{168821} i, \frac{210 \sqrt{337642}}{168821}\right)\right)$. How the different tetrahedra in the resulting finite orbit (of length 25)  map onto each other is illustrated in Figure~\ref{fig:graph_Sommerville_mod}. The orbit of this tetrahedron is displayed in  Figure~\ref{fig:Orbit_Sommerville_mod}. 
\end{example}

Consider the perturbation $\tilde\tau^S$ of the Sommerville tetrahedron defined in Example \ref{ex:Som-mod}
. This size of this perturbation is $d_{\mathcal T}(\tau^S,\tilde\tau^S)=7.75\cdot 10^{-2}$ and its orbit is shown in Figure~\ref{fig:graph_Sommerville_mod}.
%Consider the perturbation $\tilde \tau_S$ of the Sommerville tetrahedron $\tau_S$ of the size $d_{\mathcal T}(\tau^S,\tilde\tau^S)=7.75\cdot 10^{-2}$ described in Example \ref{ex:1} and displayed in red in Figure~\ref{fig:graph_Sommerville_mod}.  
We observe that the orbit of $\tilde\tau^S$ gets attracted in a sub-orbit of length 24 composed of three clusters. Similarly to what happens for $\tau^S$, the orbit of $\tilde\tau^S$ gets attracted into three clusters corresponding to the other elements of the orbit of $\tau^S$. Using the notation introduced in Proposition \ref{prop:Sommerville} and discussed in Remark \ref{rem:names}, the three clusters are referred to as $c^H = \{\tau_2,\dots,\tau_{11}\}$, $c^P = \{\tau_{12},\dots,\tau_{17}\}$ , and $c^Q = \{\tau_{18},\cdots,\tau_{25}\}$, where the elements are defined in Figure \ref{fig:graph_Sommerville_mod}. For each of these clusters, we can compute the average distance of the elements to reference shapes defining the clusters. For example, if we consider $c^H$, the average distance to $\tau^H$ is given by 
\[
\overline{d_{\mathcal T}}(c^H,\tau^H):=\frac{1}{|c^H|}\sum_{\tau\in c^H} d_{\mathcal T}(\tau,\tau^H)
\]
and similarly for other clusters. We obtain $\overline{d_{\mathcal T}}(c^H,\tau^H) = 8.41\cdot 10^{-2}$, $\overline{d_{\mathcal T}}(c^P,\tau^P) = 8.00\cdot 10^{-2}$ and $\overline{d_{\mathcal T}}(c^Q,\tau^Q) = 7.43\cdot 10^{-2}$, which are of similar size to the original perturbation.
%In Figure~\ref{fig:graph_Sommerville_mod}, we show how the tetrahedra  map onto each other under LEB. The red node in the graph corresponds to the perturbation of size $d_{\mathcal T} = 7.75\cdot 10^{-2}$.  of the Sommerville tetrahedron and does not have any incoming edges, showing that the orbit never returns to this tetrahedron. 
\begin{remark}
    In Example \ref{ex:Som-mod}, we obtained orbits of lengths $24$ and $25$. These lengths were not reported in Padrón \emph{et al.} \cite{PADRON2025555}, where they identified finite orbits of length $n\in\{4,8,9,13,21,37\}$ and the orbit of length $4$ is that of $\tau^S$. However, the normalized tetrahedron $\tilde \tau^S \in R_1^+$ (the  $R_1^+$ family is defined in Definition 2 of \cite{PADRON2025555} ), the family of tetrahedra to which their result applies. We believe that the discrepancy originates in the different normalization procedures used in both approaches discussed in Remark \ref{rem:norm}.
\end{remark}

To test whether shorter orbits are possible as well, we modified the element $\tau^H=\left(\frac12+\frac{1}2i,\left(\frac12+\frac{1}2i,\frac12\right)\right)$  from the orbit of the Sommerville tetrahedron. 

\begin{conjecture}\label{conj:1}
    Let $\tau^H = (z^H,(w^H,t^H) = \left(\frac12+\frac{1}2i,\left(\frac12+\frac{1}2i,\frac12\right)\right)$. When modifying $\tau^H$ component-wise, we can obtain orbits of lengths $n\in \{7,8,21,22,23,24,25,26\}$. More precisely, we have
    \begin{itemize}
        \item For $\tilde\tau^H_{1,\alpha} = \left(\frac12-\alpha+\frac{1}2i,\left(\frac12+\frac{1}2i,\frac12\right)\right)$, $\alpha\in ]0,\frac{1}{10}]$, the orbit $\Gamma_{\tilde\tau^H_{1,\alpha} }$  has length 21.
        \item For $\tilde\tau^H_{2,\alpha} = \left(\frac12+(\frac{1}2+\alpha)i,\left(\frac12+\frac{1}2i,\frac12\right)\right)$, $\alpha\in [-\frac{3}{50},\frac{1}{10}]$, $\alpha\neq 0$, the orbit $\Gamma_{\tilde\tau^H_{2,\alpha} }$  has length between 22 or 24.
        \item For $\tilde\tau^H_{3,\alpha} = \left(\frac12+\frac{1}2i,\left(\frac12-\alpha+\frac{1}2i,\frac12\right)\right)$, $\alpha\in ]0,\frac{2}{10}]$, the orbit $\Gamma_{\tilde\tau^H_{3,\alpha} }$ has length 23, 24 or 26.
        \item For $\tilde\tau^H_{4,\alpha} = \left(\frac12+\frac{1}2i,\left(\frac12+(\frac{1}2+\alpha)i,\frac12\right)\right)$, $\alpha\in [-\frac{2}{10},\frac{2}{10}]$, $\alpha\neq 0$, the orbit $\Gamma_{\tilde\tau^H_{4,\alpha} }$  has length 23, 24, 25, or 26.
        \item For $\tilde\tau^H_{5,\alpha} = \left(\frac12+\frac{1}2i,\left(\frac12+\frac{1}2i,\frac12+\alpha\right)\right)$, $\alpha\in [-\frac{9}{100},\frac{2}{10}]$, $\alpha\neq 0$, the orbit $\Gamma_{\tilde\tau^H_{5,\alpha} }$  has length 7 or 8.
    \end{itemize}
\end{conjecture}
%\begin{remark}
%    The maximal distances to $\tau^H$ for the different cases in Conjecture \ref{conj:1} are not constant. They vary between $d_{\mathcal T}(\tilde \tau_{2,-\frac{3}{50}},\tau^H) = 1.28\cdot 10^{-1}$ and $d_{\mathcal T}(\tilde \tau_{3,\frac{2}{10}},\tau^H) = 3.97\cdot 10^{-1}$. This suggests that the metric introduced here may not be the best one to construct a proof in the spirit of Perdomo and Plaza \cite{Perdomo2013,Perdomo2014}. Their approach relies on the fact that a union of hyperbolic circles is closed under LEB, which does not seem to be the case with the distance chosen here.
%\end{remark}

We have not been able to prove Conjecture \ref{conj:1}, but we have strong evidence that it holds true. To show that we performed experiments in which the value of $\alpha$ was systematically changed by steps of $\frac{1}{100}$ (we use $\dots$, when all the $\frac{1}{100}$ steps between limits lead to the same length of the orbit) and recorded the length of the orbit. We obtained the following results:
\begin{itemize}
    \item For $\alpha \in\{\frac{1}{100},\dots,\frac{10}{100}\}$, the orbit $\Gamma_{\tilde\tau^H_{1,\alpha}}$ has length 21 and the orbit $\Gamma_{\tilde\tau^H_{2,\alpha}}$ has length 24.
    %\item For $\alpha \in\{\frac{1}{100},\dots,\frac{10}{100}\}$, the orbit $\Gamma_{\tilde\tau^H_{2,\alpha}}$ has length 24.
    \item For $\alpha \in\{-\frac{1}{100},\dots,-\frac{3}{100}\}$, the orbit $\Gamma_{\tilde\tau^H_{2,\alpha}}$ has length 22.
    \item For $\alpha \in\{-\frac{4}{100},\dots,-\frac{6}{100}\}$, the orbit $\Gamma_{\tilde\tau^H_{2,\alpha}}$ has length 23.
    \item For $\alpha \in\{\frac{1}{100},\frac{3}{100},\frac{5}{100},\dots,\frac{9}{100}, \frac{11}{100},\frac{12}{100}, \frac{15}{100},\dots,\frac{20}{100}\}$, the orbits $\Gamma_{\tilde\tau^H_{3,\alpha}}$ and $\Gamma_{\tilde\tau^H_{4,\alpha}}$ have length 26.
    \item For $\alpha \in\{\frac{2}{100},\frac{4}{100},\frac{14}{100}\}$, the orbits $\Gamma_{\tilde\tau^H_{3,\alpha}}$ and $\Gamma_{\tilde\tau^H_{4,\alpha}}$ have length 23.
    \item For $\alpha \in\{\frac{10}{100},\frac{13}{100}\}$, the orbits $\Gamma_{\tilde\tau^H_{3,\alpha}}$ and $\Gamma_{\tilde\tau^H_{4,\alpha}}$ have length 24.
    %\item For $\alpha \in\{\frac{1}{100},\frac{3}{100},\frac{5}{100},\dots,\frac{9}{100},\frac{11}{100},\frac{12}{100},\frac{15}{100},\dots,\frac{20}{100}\}$, the orbit $\Gamma_{\tilde\tau^H_{4,\alpha}}$ has length 26.
    %\item For $\alpha \in\{\frac{2}{100},\frac{4}{100},\frac{14}{100}\}$, the orbit $\Gamma_{\tilde\tau^H_{4,\alpha}}$ has length 23.
    %\item For $\alpha \in\{\frac{10}{100},\frac{13}{100}\}$, the orbit $\Gamma_{\tilde\tau^H_{4,\alpha}}$ has length 24.
    \item For $\alpha \in\{-\frac{1}{100},-\frac{3}{100},-\frac{5}{100},\dots,-\frac{7}{100},-\frac{9}{100},-\frac{11}{100},-\frac{12}{100},-\frac{15}{100},\dots,\frac{19}{100}\}$, the orbit $\Gamma_{\tilde\tau^H_{4,\alpha}}$ has length 23.
    \item For $\alpha =\frac{2}{10}$, the orbit $\Gamma_{\tilde\tau^H_{4,\alpha}}$ has length 25.
    \item For $\alpha \in\{\frac{1}{100},\dots,\frac{20}{100}\}$, the orbit $\Gamma_{\tilde\tau^H_{5,\alpha}}$ has length 8.
    \item For $\alpha \in\{-\frac{1}{100},-\frac{3}{100},-\frac{5}{100},\dots,-\frac{9}{100}\}$, the orbit $\Gamma_{\tilde\tau^H_{5,\alpha}}$ has length 8.
    \item For $\alpha \in\{-\frac{2}{100},-\frac{34}{100}\}$, the orbit $\Gamma_{\tilde\tau^H_{5,\alpha}}$ has length 7.
\end{itemize}

Outside of the intervals specified in Conjecture \ref{conj:1}, the orbits seem to diverge. For example, we observe that the length of $\Gamma_{\tilde\tau^H_{2,\frac{11}{100}}}$ already contains 382 shapes after 20 iterations. 

\begin{figure}
    \centering
\begin{tabularx}{0.48\textwidth}{@{}lX@{}}
\hline
ID & Normalization \\
\hline
$\tau_{1}$ &  \begin{minipage}{\linewidth}$\left(\frac{1}{2} + \frac{11}{20} i, \left(\frac{1}{2} + \frac{11}{20} i, \frac{1}{2}\right)\right)$\end{minipage} \\
$\tau_{2}$ &  \begin{minipage}{\linewidth}$\left(\frac{1}{2} + \frac{1}{2} i, \left(\frac{1}{2} + \frac{1}{2} i, \frac{11}{20}\right)\right)$\end{minipage} \\
$\tau_{3}$ &  \begin{minipage}{\linewidth}$\left(\frac{1}{2} + \frac{5}{11} i, \left(\frac{1}{2} + \frac{5}{11} i, \frac{5}{11}\right)\right)$\end{minipage} \\
$\tau_{4}$ &  \begin{minipage}{\linewidth}$\left(\frac{100}{321} + \frac{10 \sqrt{221}}{321} i, \left(\frac{200}{321} + \frac{1210 \sqrt{221}}{70941} i, \frac{110 \sqrt{70941}}{70941}\right)\right)$\end{minipage} \\
$\tau_{5}$ &  \begin{minipage}{\linewidth}$\left(\frac{100}{321} + \frac{10 \sqrt{221}}{321} i, \left(\frac{221}{321} + \frac{1000 \sqrt{221}}{70941} i, \frac{110 \sqrt{70941}}{70941}\right)\right)$\end{minipage} \\
$\tau_{6}$ &  \begin{minipage}{\linewidth}$\left(\frac{1}{2} + \frac{11 \sqrt{2}}{40} i, \left(\frac{1}{2} + 0 i, \frac{1}{2}\right)\right)$\end{minipage} \\

$\tau_{7}$ &  \begin{minipage}{\linewidth}$\left(\frac{1}{2} + \frac{5 \sqrt{221}}{221} i, \left(\frac{100}{221} + 0 i, \frac{110}{221}\right)\right)$\end{minipage} \\
$\tau_{8}$ &  \begin{minipage}{\linewidth}$\left(\frac{1}{2} + \frac{5 \sqrt{221}}{221} i, \left(\frac{121}{221} + 0 i, \frac{110}{221}\right)\right)$\end{minipage} \\
\hline
\end{tabularx}
\raisebox{-.5\height}{
    \includegraphics[width=0.45\linewidth]{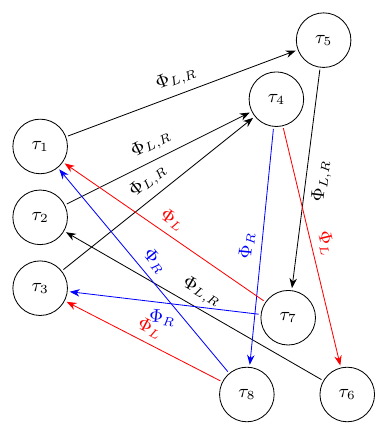}}
    \caption{An example of an orbit of length 8. Left: list of tetrahedra in the orbit. Right: graph of how they map onto each other.}
    \label{fig:table_and_graph8}
\end{figure}

\begin{figure}
    \centering
    \includegraphics[width=0.85\linewidth]{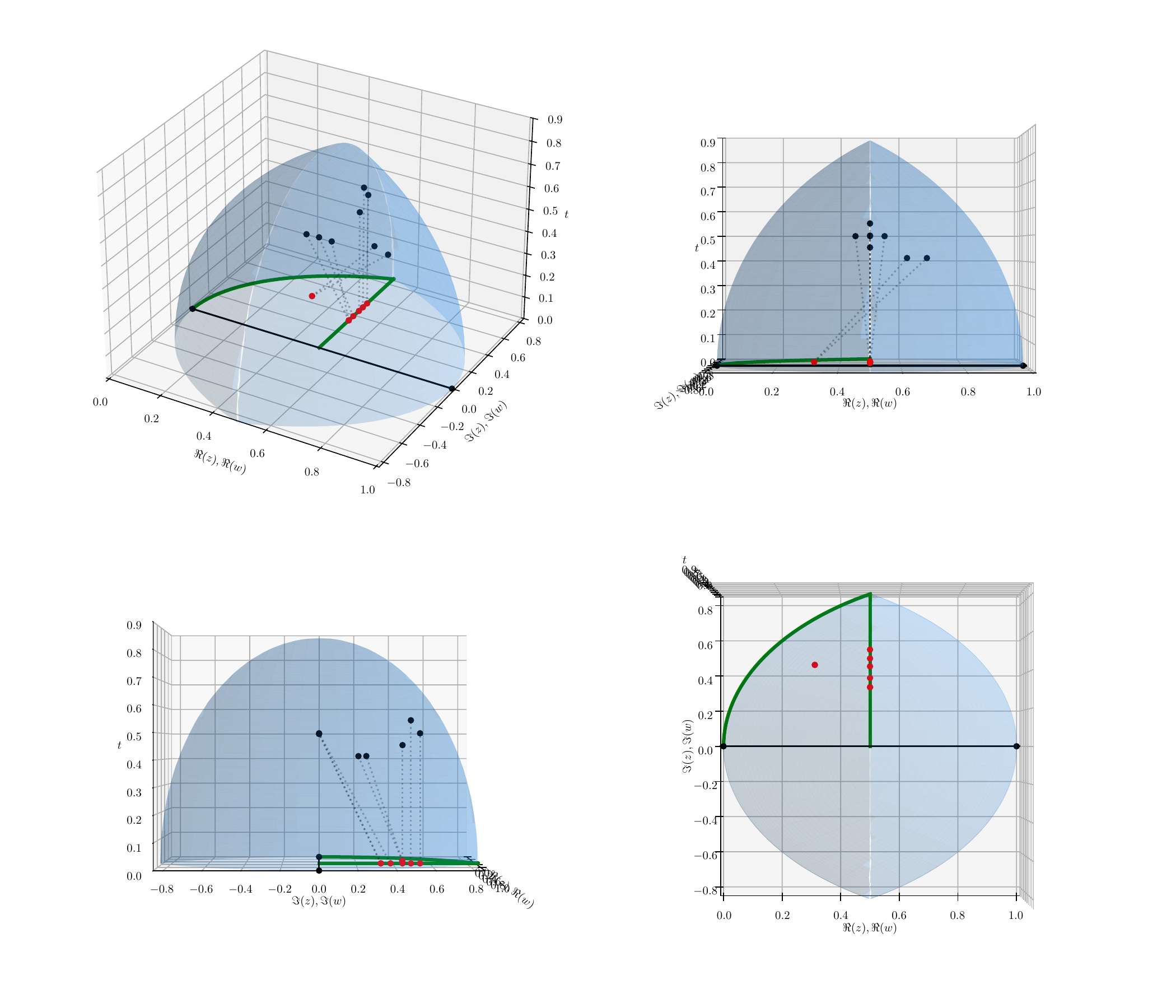}
    \caption{Visualization of the orbit of length 8 from Figure~\ref{fig:table_and_graph8}. The bottom right panel shows only the positions of the base vertices.}
    \label{fig:viz_orbit_8}
\end{figure}

\begin{example}\label{ex:orbit8}
    We modify the tetrahedron with normalization $\tau^H =\left(\frac12+\frac{1}2i,\left(\frac12+\frac{1}2i,\frac12\right)\right)$ from the orbit of the Sommerville tetrahedron  to $\tilde\tau^H_{5,\frac{1}{20}}=
    \left(\frac{1}{2} + \frac{1}{2} i, \left(\frac{1}{2} + \frac{1}{2} i, \frac{11}{20}\right)\right)$. This tetrahedron has an orbit of length 8 containing the tetrahedra listed in Figure~\ref{fig:table_and_graph8} and visualized in Figure~\ref{fig:viz_orbit_8}.
\end{example}

Example \ref{ex:orbit8} illustrates the orbit of $\tilde \tau^H_{5,\frac{1}{20}}$, which has an orbit of length 8 as predicted by Conjecture \ref{conj:1}. This is a perturbation of size $d_{\mathcal{T}}(\tilde \tau^H_{5,\frac{1}{20}},\tau^H)=9.53\cdot 10^{-2}$. In Figure~\ref{fig:table_and_graph8}, the 8 shapes in the orbit of $\tilde \tau^H_{5,\frac{1}{20}}$ are listed in the table on the left, and the graph of how they map onto each other is displayed on the right. A visualization of the orbit is provided in Figure \ref{fig:viz_orbit_8}. We clearly identify three clusters of normalized tetrahedra corresponding to $\tau^H$, $\tau^P$ and $\tau^Q$ denoted as before by $c^H_{5,\frac{1}{20}} = \{\tau_1,\tau_2,\tau_3\}$, $c^P_{5,\frac{1}{20}} = \{\tau_4,\tau_5\}$, and $c^Q_{5,\frac{1}{20}} = \{\tau_6,\tau_7,\tau_8\}$. The average distances to the reference elements are $\overline{d_{\mathcal T}}(c^H_{5,\frac{1}{20}},\tau^H) = 1.33\cdot 10^{-1}$, $\overline{d_{\mathcal T}}(c^P_{5,\frac{1}{20}},\tau^P) = 1.12\cdot 10^{-1}$, and $\overline{d_{\mathcal T}}(c^Q_{5,\frac{1}{20}},\tau^Q) = 1.03\cdot 10^{-1}$.

%A short exploration yielded orbits of lengths $n\in\{3,4,8, 13, 14, 24, 25,  26\}$.  These lengths are overlapping with those obtained by Padrón \emph{et al.}\cite{PADRON2025555}, but not identical. Since we already have an example of an orbit of length 24, we will show an example of an orbit of length 8. Example~\ref{ex:orbit8} provides such an example. The results are shown in Figures~\ref{fig:table_and_graph8} and \ref{fig:viz_orbit_8}. 

\begin{remark}\label{rem:orbit-lenght}
    Conjecture \ref{conj:1} only concerns modifications of a single component of $\tau^H$. When modifying more components at once, we can obtain orbits of other lengths. For example, the tetrahedron $\tilde\tau_1=\left(\frac{1}{2} + \frac{\sqrt{442}}{40} i, \left(\frac{1}{2} + \frac{5 \sqrt{442}}{221} i, \frac{11 \sqrt{442}}{442}\right)\right)$ has an orbit orbit of length 13 that can be reached from the modification of the Sommerville tetrahedron $\tilde\tau^S_1 = \left(\frac{100}{221} + \frac{5 \sqrt{926}}{221} i, \left(\frac{100}{221} + \frac{105 \sqrt{926}}{102323} i, \frac{220 \sqrt{102323}}{102323}\right)\right)$ which has an orbit of length 14. Furthermore, the tetrahedron $\tilde\tau_2=\left(\frac{1}{2} + \frac{1}{2} i, \left(\frac{53}{100} + \frac{49}{100} i, \frac{3}{5}\right)\right)$ has an orbit of length 44.
\end{remark}

The results obtained in this section show that small perturbations of $\tau^S$ or of $\tau^H$, may lead to finite orbits of various lengths. So far, we have shown orbits of lengths $n\in\{3,4,7,8,13,14,21,22,23,24,25,26,44\}$. In the next section, we look at small perturbations of the regular tetrahedron, since this family of tetrahedra has been extensively studied in \cite{Suarez2021,Trujillo2024}.
%These two [representative ??] examples suggest that the finite sub-orbit of the orbit of the Sommerville tetrahedron plays a central role in understanding the dynamics of LEB on [most of ???] tetrahedron. [However, larger ??] perturbations of the Sommerville tetrahedron do not seem to generate finite orbits. In the next section, we study tetrahedra that do not seem to have  [finite orbitS ??]. 

\subsection{Orbits of small perturbations of the regular tetrahedron}
In \cite{Suarez2021,Trujillo2024,PADRON2025555}, a family of nearly equilateral tetrahedra is studied and it is shown that at most 37 similarity classes can be produced. This limit was mentioned by Adler \cite{adler1983bisection} already in 1983, but without any proof. Here we apply our normalization procedure to six tetrahedra listed in \cite{Suarez2021} to test whether our normalization procedure and our implementation the LEB algorithm reproduces the results mentioned above. Note that the orbit of length 44 mentioned in Remark \ref{rem:orbit-lenght} does not contradict the upper limit of 37 similarity classes, since the tetrahedron taken as an example does not belong to the family of tetrahedra studied in \cite{Suarez2021,Trujillo2024,PADRON2025555}. 

\begin{table}[ht]
\small
\centering
\caption{List of the tetrahedra used in \cite{Suarez2021} and the length of their orbits.}
\begin{tabularx}{\textwidth}{l p{.35\textwidth} p{.38\textwidth} c}
\toprule
\bf ID &\bf  Vertex coordinates & \bf Normalized form & \bf Orbit length \\
\midrule
T1 &
$\begin{aligned}\,&(0,0,0),\;(\sqrt{22},0,0),\;\big(\tfrac{5\sqrt{22}}{11},\tfrac{\sqrt{1507}}{11},0\big),\\
\;&\big(\tfrac{21\sqrt{22}}{44},\tfrac{49\sqrt{1507}}{3014},\tfrac{\sqrt{3480622}}{548}\big)\,\end{aligned}$
& $\left(\frac{5}{11} + \frac{\sqrt{274}}{22} i, \left(\frac{21}{44} + \frac{49 \sqrt{274}}{6028} i, \frac{\sqrt{19143421}}{6028}\right)\right)$&39\\[20pt]

T2 &
$\begin{aligned}&\,(0,0,0),\;(\sqrt{21},0,0),\;\big(\tfrac{3\sqrt{21}}{7},\tfrac{\sqrt{595}}{7},0\big),\\&
\;\big(\tfrac{10\sqrt{21}}{21},\tfrac{31\sqrt{595}}{1190},\tfrac{\sqrt{3078105}}{510}\big)\,\end{aligned}$
& $\left(\frac{3}{7} + \frac{\sqrt{255}}{21} i, \left(\frac{10}{21} + \frac{31 \sqrt{255}}{3570} i, \frac{\sqrt{7182245}}{3570}\right)\right)$&37\\[20pt]

T3 &
$\begin{aligned}&\,(0,0,0),\;(\sqrt{21},0,0),\;\big(\tfrac{\sqrt{21}}{2},\tfrac{\sqrt{47}}{2},0\big),\\&
\;\big(\tfrac{11\sqrt{21}}{21},\tfrac{4\sqrt{47}}{47},\tfrac{2\sqrt{2654043}}{987}\big)\end{aligned}$
& $\left(\frac{10}{21} + \frac{2 \sqrt{59}}{21} i, \left(\frac{1}{2} + \frac{\sqrt{59}}{59} i, \frac{\sqrt{3331671}}{2478}\right)\right)$&43 \\[20pt]

T4 &
$\begin{aligned}&(0,0,0),\;(\sqrt{41},0,0),\big(\tfrac{20\sqrt{41}}{41},\tfrac{\sqrt{47478}}{41},0\big),\\&
\;\big(\tfrac{21\sqrt{41}}{41},\tfrac{53\sqrt{47478}}{7913},\tfrac{\sqrt{935471}}{193}\big)\,\end{aligned}$
& $\left(\frac{20}{41} + \frac{\sqrt{1117}}{41} i, \left(\frac{21}{41} + \frac{318 \sqrt{1117}}{45797} i, \frac{\sqrt{1331868354}}{45797}\right)\right)$&43 \\[20pt]

T5 &
$\begin{aligned}&(0,0,0),\;(10,0,0),\big(\tfrac{99}{20},\tfrac{\sqrt{28999}}{20},0\big),\\&
\;\big(\tfrac{101}{20},\tfrac{9001\sqrt{28999}}{579980},\tfrac{\sqrt{1870811}}{28999}\big)\end{aligned}$
& $\left(\frac{99}{200} + \frac{\sqrt{28599}}{200} i, \left(\frac{101}{200} + \frac{9001 \sqrt{28599}}{5719800} i, \frac{\sqrt{53503323789}}{285990}\right)\right)$&43 \\[20pt]

T6 &
$\begin{aligned}&(0,0,0),\;(10,0,0),\big(\tfrac{99}{20},\tfrac{\sqrt{28999}}{20},0\big),\\&
\;\big(\tfrac{101}{20},\tfrac{8801\sqrt{28999}}{579980},\tfrac{\sqrt{1850713}}{28999}\big)\end{aligned}$
& $\left(\frac{99}{200} + \frac{\sqrt{28199}}{200} i, \left(\frac{101}{200} + \frac{8801 \sqrt{28199}}{5639800} i, \frac{\sqrt{52188255887}}{281990}\right)\right)$&37 \\

\bottomrule
\end{tabularx}
\label{tab:T1-6}
\end{table}

\begin{theorem}\label{thm:nearly-equi}
    The six tetrahedra listed in Table \ref{tab:T1-6} have finite orbits as listed in the table.
\end{theorem}
\begin{proof}
    This is shown by direct computation. For example, the normalized forms of the tetrahedra in the orbit of T1 are listed in Table~\ref{tab:tet_nodes_39}.
\end{proof}
\begin{remark}
    Theorem \ref{thm:nearly-equi} shows that, contrarily to the conjecture by Adler \cite{adler1983bisection} and studied in \cite{Trujillo2024,PADRON2025555}, it is possible to obtain finite orbits longer than 37. We do recover, for T2 and T6, the orbit of length 37, but the other tetrahedra in Table~\ref{tab:T1-6} have longer orbits. We believe that the discrepancy comes from the different normalization used in our work from that one used in \cite{Trujillo2024,PADRON2025555}, which, as pointed out in Remark \ref{rem:norm} may lead to a different choices of the longest edge to be bisected and hence, may lead to orbits of different lengths. As a result, the proofs developed in \cite{Trujillo2024,PADRON2025555} are only valid for their normalization procedure and the implied choice of the edge to be bisected. In general, their result that the number of similarity classes of nearly equilateral tetrahedra being bounded by $37$ does not seem to hold.
\end{remark}

\begin{figure}
    \centering
    \includegraphics[width=0.85\linewidth]{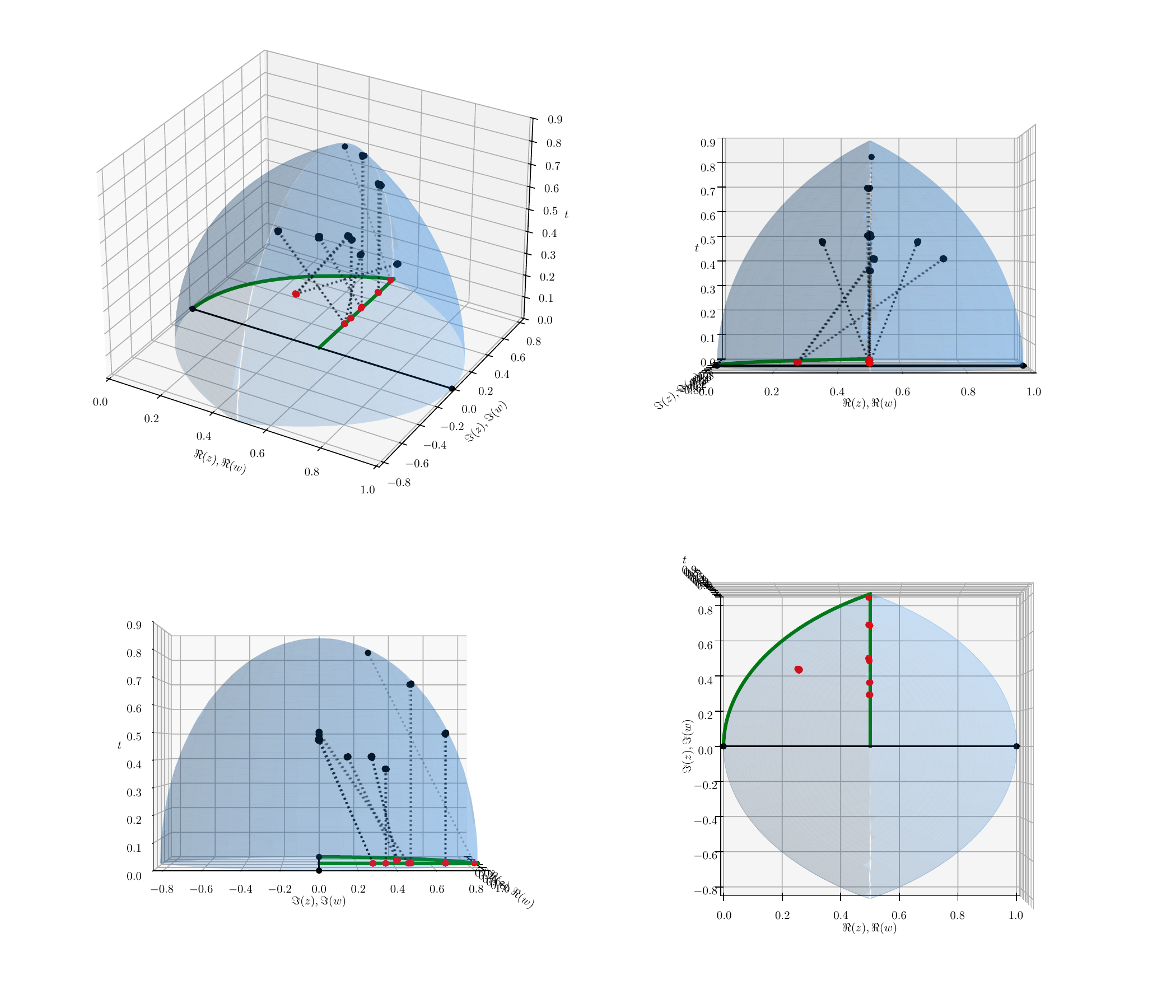}
    \caption{Orbit of the tetrahedron T5 from Table~\ref{tab:T1-6}. The bottom right panel shows only the positions of the base vertices.}
    \label{fig:orbitT5}
\end{figure}

An interesting question is how do these orbits look like, do they also cluster like the perturbations of $\tau^S$ and $\tau^H$? To answer this question, we provide a visualization of the orbit of T5 containing 43 shapes in Figure \ref{fig:orbitT5}.

As can be seen from the figure, the orbit clearly shows clustering. The 43 shapes seem to be attracted to a sub-orbit of length 42 that clusters around 8 points in $\mathcal T$. This number of clusters is reminiscent of the result in \cite{aparicio2015minimum} which states that the minimal number of similarity classes of the regular tetrahedron following the Branch-and-Bound algorithm is 8. Since the tetrahedron is nearly equilateral, some of the ambiguity of choosing the edge is eliminated, thus allowing the orbit to be attracted to a short orbit. 
\begin{remark}\label{rem:orbit_lengths2}
The fact that the orbit never comes back to its starting point is valid for all tetrahedra listed in Table \ref{tab:T1-6}. This implies the existence of orbits of length $n\in \{36,37,38,39,42,43\}$. 
\end{remark}
\begin{remark}
    The attracting sub-orbit in the orbit of T5 does not overlap that of $\tau^P$. These are two different attracting cycles. 
\end{remark}
\begin{theorem}\label{thm:length}
    There exists finite orbits under LEB of lengths $n\in\{3,4,7,8,13,14\} \cup \{21,\dots,26\}\cup\{36,37,38,39,42,43,44\}$.
\end{theorem}
\begin{proof}
    This Theorem collects the results of Thm \ref{prop:Sommerville}, Examples \ref{ex:Som-mod} and \ref{ex:orbit8}, of the illustration of Conjecture \ref{conj:1} and Remaks \ref{rem:orbit-lenght} and \ref{rem:orbit_lengths2}.
\end{proof}
In the next section, we will show that the regular tetrahedron is not as well-behaved and that its orbit does not seem to be finite.

\subsection{On the behaviour of LEB on general tetrahedra}
\begin{figure}[t]
    \centering
    \includegraphics[width=0.65\linewidth]{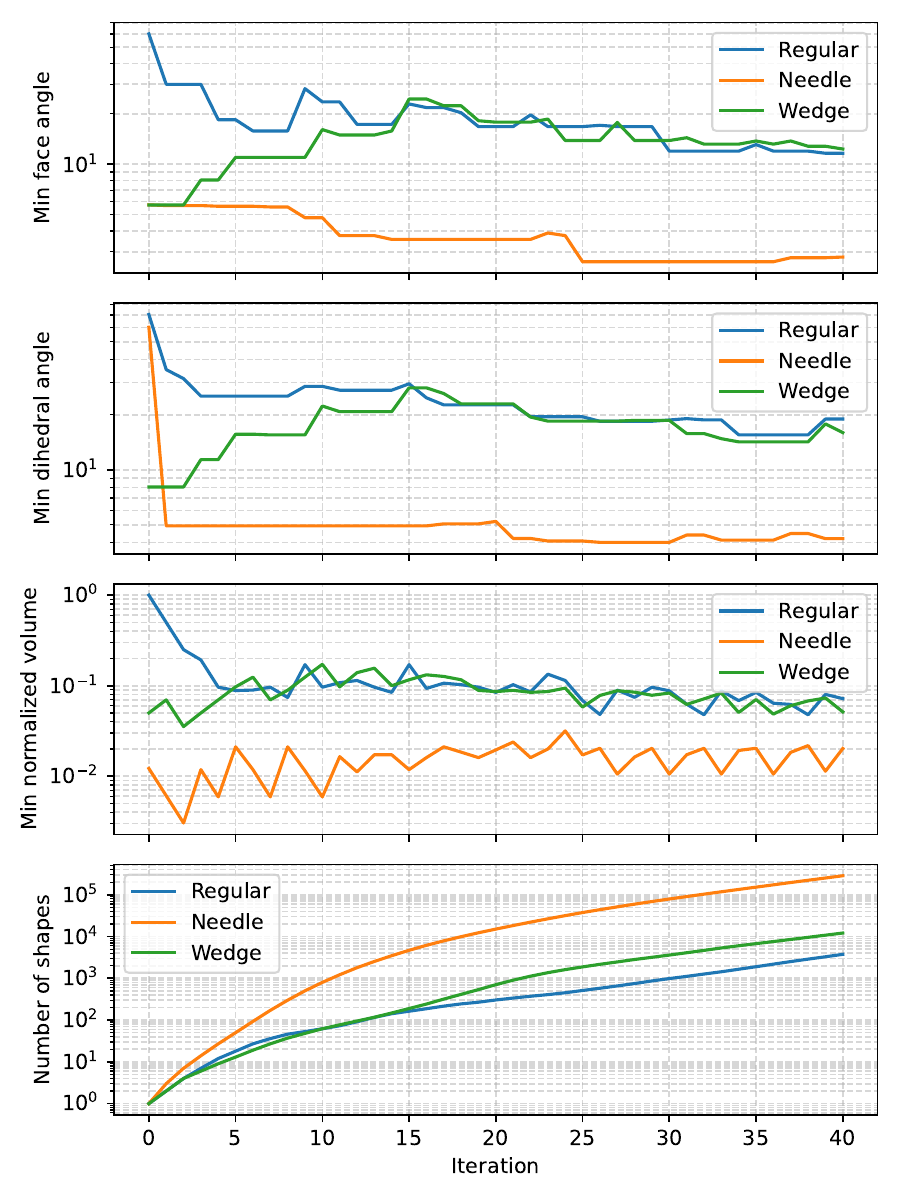}
    \caption{Dynamics of the metrics during the LEB process. The metrics at each iteration are computed on the set of new shapes produced at that iteration.}
    \label{fig:dynamics}
\end{figure}

We now discuss four cases of tetrahedra that do not seem to have finite orbits. We consider the regular tetrahedron, the cube corner tetrahedron and two nearly degenerate tetrahedra and iterate the LEB process for $n=40$ iterations and record how many shapes appear in the resulting orbits, the minimal face and dihedral angles as well as their normalized volumes. In order to get a better interpretation of the volumes, we normalize the volume by the maximum possible volume of a normalized tetrahedron, i.e. the volume of the regular tetrahedron of edge length 1. This volume is given by $\frac{\sqrt{2}}{12}$. Hence, the volumes are reported in percent of the volume of the regular tetrahedron, see Table \ref{tab:leb-tetra} and Figure \ref{fig:dynamics}.

\begin{example}\label{ex:regular}
    The regular tetrahedron is given in normalized form as $\tau^R=\left(\frac{1}{2} + \frac{\sqrt{3}}{2} i, \left(\frac{1}{2} + \frac{\sqrt{3}}{6} i, \frac{\sqrt{6}}{3}\right)\right)$. Applying the LEB refinement process to this tetrahedron does not seem to lead to a finite orbit. 
\end{example}

\begin{example}\label{ex:cubecorner}
    The cube corner tetrahedron with vertices $(0,0,0), (1,0,0), (0,1,0)$ and $(0,0,1)$ is given in normalized form as $\tau^C=\left(\frac{1}{2} + \frac{1}{2} i, \left(\frac{1}{2} + \frac{1}{2} i, \frac{\sqrt{2}}{2}\right)\right)$. Applying the LEB refinement process to this tetrahedron does not seem to lead to a finite orbit. 
\end{example}

\begin{example}\label{ex:wedge}
    The wedge tetrahedron with vertices $(0,0,0), (1,0,0), (0,1,0)$ and $(0,0,\frac{1}{10})$ is given in normalized form as $\tau^W=\left(\frac{1}{2} + \frac{1}{2} i, \left(\frac{1}{2} + \frac{1}{2} i, \frac{\sqrt{2}}{20}\right)\right)$. Applying the LEB refinement process to this tetrahedron does not seem to lead to a finite orbit. 
\end{example}

\begin{example}\label{ex:needle}
    The needle tetrahedron also with vertices $(0,\frac{1}{10\sqrt{3}},0), (\frac{1}{20},-\frac{1}{20\sqrt{3}},0),(-\frac{1}{20},-\frac{1}{20\sqrt{3}},0)$ and $(0,0,1)$ (its base is an equilateral triangle in $xy-$plane with the edge lengths equal to 0.1) used in \cite{Hannukainen2014} is given in normalized form as $\tau^N=\left(\frac{3}{602} + \frac{\sqrt{3603}}{602} i, \left(\frac{3}{602} + \frac{599 \sqrt{3603}}{723002} i, \frac{30 \sqrt{1084503}}{361501}\right)\right)$. Applying the LEB refinement process to this tetrahedron does not seem to lead to a finite orbit.
\end{example}

The four tetrahedra given in Examples~\ref{ex:regular}, \ref{ex:cubecorner}, \ref{ex:wedge}, and \ref{ex:needle} have long and potentially unbounded orbits. The regular tetrahedron $\tau^R$ and cube corner tetrahedron $\tau^C$ have large orbits, but these orbits do not grow uncontrollably as shown in Table~\ref{tab:leb-tetra} since after $n=40$ the orbits $\bigcup_{n=0}^{40}\Gamma^{(n)}_{\tau^R}$ and $\bigcup_{n=0}^{40}\Gamma^{(n)}_{\tau^C}$ contain $3793$ and $2919$ tetrahedra, respectively. In addition, the minimal face and dihedral angles seem to remain bounded by similar values. 

% --- table ---
\begin{table}[t]
  \centering
  \caption{We report the minimal dihedral angle, the minimal face angle and the minimal normalized volume of the tetrahedra in the orbits $\bigcup_{n=0}^{40}\Gamma^{(n)}_{\tau}$, where $\tau$ is one of the normalized tetrahedra from Examples~\ref{ex:regular}, \ref{ex:cubecorner}, \ref{ex:wedge} and \ref{ex:needle}. The initial values are given in parentheses.}
  \label{tab:leb-tetra}
  % Column layout:
  %  - first column: left-aligned flexible X (metric name)
  %  - other three columns: centered flexible X (one per tetrahedron case)
  \begin{tabularx}{\textwidth}{
      >{\raggedright\arraybackslash}X
      >{\centering\arraybackslash}X
      >{\centering\arraybackslash}X
      >{\centering\arraybackslash}X
      >{\centering\arraybackslash}X
    }
    \toprule
    Metric & Regular    & Cube corner & Wedge& Needle \\
    \midrule
    Minimal dihedral angle (deg)        & \(15.52^\circ \ (70.53^\circ) \)  & \(15.52^\circ\ (54.74^\circ)\) &\(8.05^\circ\ (8.05^\circ)\)& \(4.01^\circ\ (60.08^\circ)\)\\
    Minimal face angle (deg)            & \(11.63^\circ\ (60^\circ) \)  & \(11.97^\circ\ (45^\circ)\)&\(5.71^\circ\ (5.71^\circ)\)& \(2.62^\circ\ (5.72^\circ) \) \\
    Minimal normalized volume (\(\%\)) & \(4.78\%\ (100\%)\)  & \(4.78\%\ (50\%)\)& \(3.54\%\ (5\%)\)& \(0.30\%\ (1.22\%)\) \\
    Number of shapes in the orbit after $n=40$ iterations       & \(3793\)  & \(2919 \) & \(12144 \)& \(285592 \)\\
    % add more rows as needed
    \bottomrule
  \end{tabularx}
\end{table}

In Theorem \ref{prop:regular-cube}, we show that the orbit of the cube corner tetrahedron is actually included in that of the regular tetrahedron, explaining why the values obtained in the Table \ref{tab:leb-tetra} are so close to each other and why minimal quantities for the cube corner are always larger than minimal quantities for the regular tetrahedron.  
\begin{theorem}\label{prop:regular-cube}
    Assuming that the edge to be bisected is specified by our normalization procedure, the orbit of the cube corner tetrahedron is included in the orbit of the regular tetrahedron.
\end{theorem}

\begin{proof}
    It suffices to show how to reach the cube corner tetrahedron $\tau_C$ from the regular tetrahedron $\tau_R$ using the maps $\Phi_L$ and $\Phi_R$. This can be done in 7 steps:
    \begin{enumerate}
        \item Step 1: $\Phi_{L,R}(\tau^R) = \left(\frac{1}{4} + \frac{\sqrt{3}}{4} i, \left(\frac{1}{2} + \frac{\sqrt{3}}{6} i, \frac{\sqrt{6}}{3}\right)\right)$
        \item Step 2: $\Phi_R\left(\left(\frac{1}{4} + \frac{\sqrt{3}}{4} i, \left(\frac{1}{2} + \frac{\sqrt{3}}{6} i, \frac{\sqrt{6}}{3}\right)\right)\right) = \left(\frac{1}{4} + \frac{\sqrt{3}}{4} i, \left(\frac{1}{2} + \frac{\sqrt{3}}{3} i, \frac{\sqrt{6}}{6}\right)\right)$
        \item Step 3: $\Phi_L\left(\left(\frac{1}{4} + \frac{\sqrt{3}}{4} i, \left(\frac{1}{2} + \frac{\sqrt{3}}{3} i, \frac{\sqrt{6}}{6}\right)\right)\right) = \left(\frac{1}{2} + \frac{\sqrt{3}}{6} i, \left(\frac{1}{3} + 0 i, \frac{\sqrt{2}}{3}\right)\right)$
        \item Step 4: $\Phi_R\left(\left(\frac{1}{2} + \frac{\sqrt{3}}{6} i, \left(\frac{1}{3} + 0 i, \frac{\sqrt{2}}{3}\right)\right)\right) = \left(\frac{1}{2} + \frac{\sqrt{2}}{4} i, \left(\frac{1}{2} + \frac{\sqrt{2}}{4} i, \frac{\sqrt{2}}{4}\right)\right)$
        \item Step 5: $\Phi_{L,R}\left(\left(\frac{1}{2} + \frac{\sqrt{2}}{4} i, \left(\frac{1}{2} + \frac{\sqrt{2}}{4} i, \frac{\sqrt{2}}{4}\right)\right)\right) = \left(\frac{1}{4} + \frac{\sqrt{3}}{4} i, \left(\frac{1}{2} + \frac{\sqrt{3}}{6} i, \frac{\sqrt{6}}{6}\right)\right)$
        \item Step 6: $\Phi_L\left(\left(\frac{1}{4} + \frac{\sqrt{3}}{4} i, \left(\frac{1}{2} + \frac{\sqrt{3}}{6} i, \frac{\sqrt{6}}{6}\right)\right)\right) = \left(\frac{1}{2} + \frac{1}{2} i, \left(\frac{1}{2} + 0 i, \frac{1}{2}\right)\right)$
        \item Step 7: $\Phi_{L,R}\left(\left(\frac{1}{2} + \frac{1}{2} i, \left(\frac{1}{2} + 0 i, \frac{1}{2}\right)\right)\right) = \left(\frac{1}{2} + \frac{1}{2} i, \left(\frac{1}{2} + \frac{1}{2} i, \frac{\sqrt{2}}{2}\right)\right) = \tau^C$.
    \end{enumerate}
    This concludes the proof, since we have shown how to reach $\tau^C$ from $\tau^R$.
    %Sequence for the proof: $\tau_{21}, L,R, \tau_3,R, \tau_{14},L,\tau_{25},R,\tau_{10},L,R,\tau_{16},L, \tau_{15},L,R,\tau_{19}$
\end{proof}

In Figure \ref{fig:dynamics}, we show how the metrics for the four tetrahedra studied here evolve under LEB. We do not report the results for the cube corner, since its orbit is included in that of the regular tetrahedron. For efficiency reasons, at each iteration we compute the different metrics on newly encountered shapes only. This is possible since we are computing the orbits symbolically. When the orbits are very large, the computational cost of computing everything symbolically is too high. In this case, we consider two tetrahedra to be the same numerically if every components of their normalized form are equal when rounded at 10 decimal places. This numerical procedure has been used to obtain the results for $\tau^W$ and $\tau^N$. The number of shapes using this numerical approach has been validated by comparing it to the symbolic computations for the regular and cube corner tetrahedra and they lead to the same metrics and to the same number of shapes.

From Figure \ref{fig:dynamics} we observe that the dynamics of the metrics for the regular tetrahedron seem to decrease and the minimum values reported in Table \ref{tab:leb-tetra} occur towards the end of the simulation. While it is possible that these values continue to decrease, they seem to stabilize. The shapes produced by the wedge tetrahedron $\tau^N$ seem to become more regular and increasingly similar to that of the shapes produced by the regular tetrahedron $\tau^R$. The fact that the minimal face and dihedral angles improve is best illustrated by noted that these minimal values are equal to their initial values, see Table \ref{tab:leb-tetra}. 

The nearly degenerate needle tetrahedron $\tau^N$ is much less well-behaved. The number of shapes in $\bigcup_{n=0}^{40}\Gamma^{(n)}_{\tau^N}$ is much larger than that of the regular or cube corner tetrahedra, but still much smaller than the theoretical maximum number of shapes in such an orbit, which is $2^{41}-1$. Even though the values for the minimal face angle and minimal dihedral angle are very small, they still seem to be bounded from below. Maybe a better sign that the shapes are not degenerating is the fact that the normalized volume of newly created shapes is increasing and the minimal value is obtained at iteration 2, see Figure \ref{fig:dynamics}. Note that while the minimal normalized volume in the orbit of $\tau^N$ is $0.30\%$ of that of the regular tetrahedron, it is only a reduction by a factor of approximately $4$ from the initial normalized volume of $\tau^N$, which is $1.22\%$.
\medskip

%The minimal volume attained by elements in the orbits of the regular and the cube corner tetrahedra is $4.78\%$ of the volume of the regular tetrahedron. In addition, the minimal volume attained by elements in the orbit of the needle tetrahedron is $0.30\%$ of the regular tetrahedron. But noting that the initial volume of the needle tetrahedron is only $1.22\%$, this is only a reduction by a factor 4.

\begin{figure}
    \centering
    \includegraphics[width=0.85\linewidth]{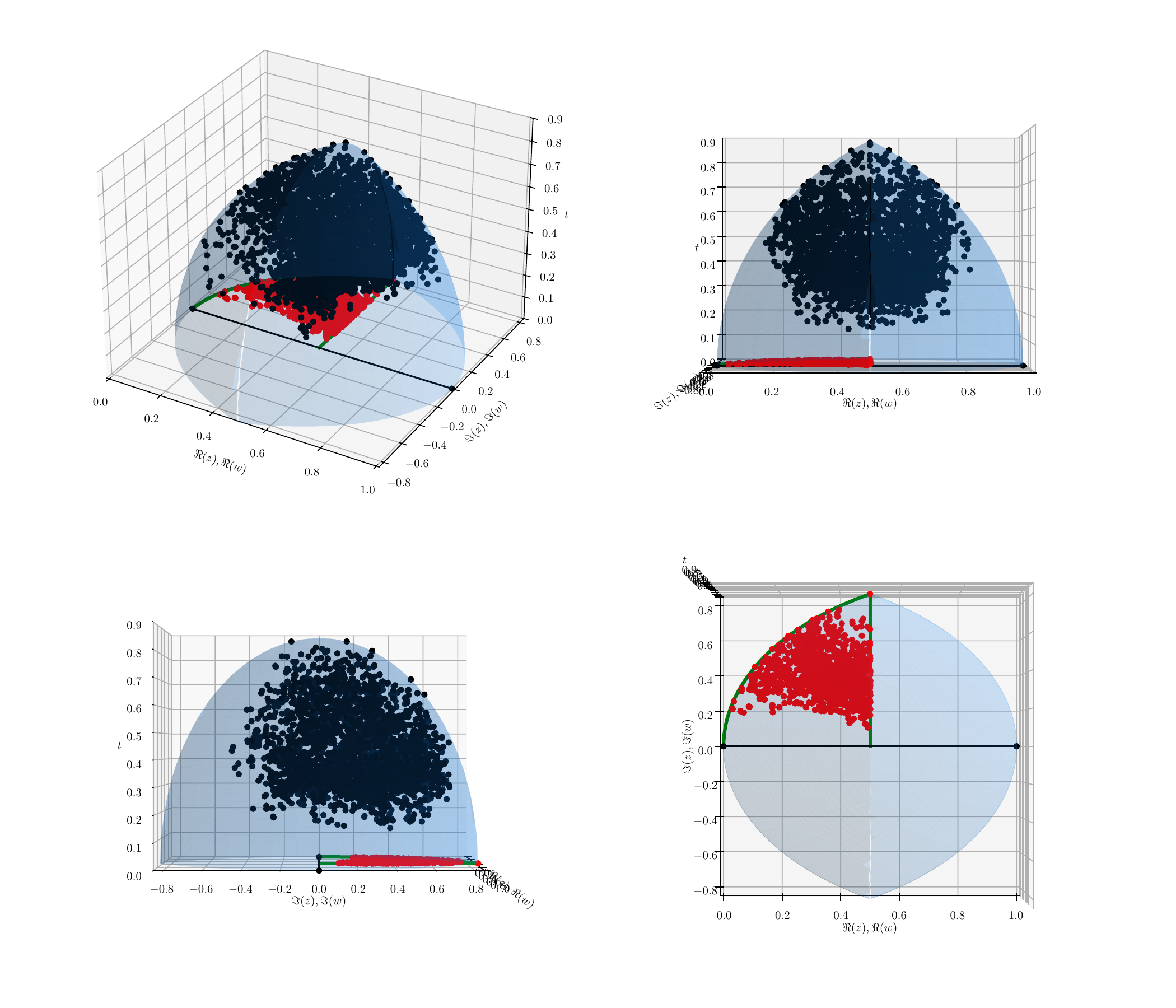}
    \caption{Visualization of the orbit of the regular tetrahedron after 40 iterations. The edges are not displayed to avoid having too much information in the figure and to better visualize the gap between the fourth vertex positions and the $xy-$plane. The bottom right panel shows only the positions of the base vertices.}
    \label{fig:orbit_regular_40}
\end{figure}

\begin{figure}
    \centering
    \includegraphics[width=0.85\linewidth]{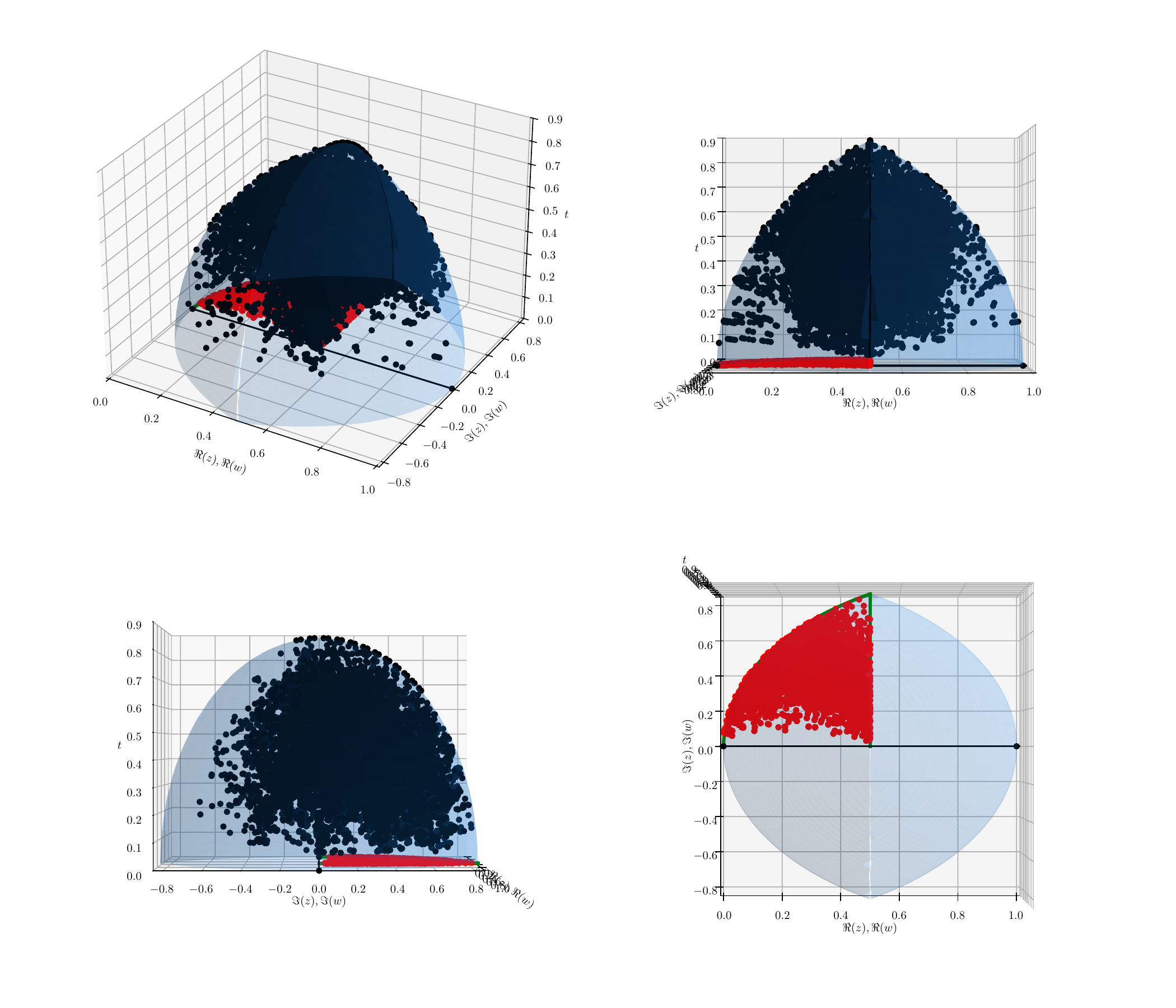}
    \caption{Visualization of the orbit of the needle tetrahedron after 20 iterations. The edges are not displayed to avoid having too much information in the figure and to better visualize the gap between the fourth vertex positions and the $xy-$plane. The bottom right panel shows only the positions of the base vertices.}
    \label{fig:orbit_needle_20}
\end{figure}

A visualization of the orbits of the regular tetrahedron and the needle tetrahedron is provided in Figures~\ref{fig:orbit_regular_40} and \ref{fig:orbit_needle_20}. A tetrahedron degenerates if either the red point corresponding to the base of the normalized tetrahedron becomes close to the $x-$axis or if the fourth vertex becomes close to the $xy-$plane. These gaps are clearly visible in Figure \ref{fig:orbit_regular_40} but much less clear in the nearly degenerate needle tetrahedron case, see Figure~\ref{fig:orbit_needle_20}.

\section{Conclusions and future work\label{sec:concl}}
In this work, we have introduced a new normalization procedure for tetrahedra inspired by the works of Perdomo and Plaza \cite{Perdomo2013,Perdomo2014}. This normalization procedure provides a geometric interpretation of the normalized form and allows for both an analytical and visual exploration of the orbits of tetrahedra under the LEB refinement algorithm. We show that the space-filling Sommerville tetrahedron essentially plays a central role in the dynamics of the LEB refinement method and that it has a very short orbit containing only 4 elements. Notably, three of these elements form a cycle. Small perturbations of the Sommerville tetrahedron lead to finite orbits clustering around this cycle. Looking at nearly regular tetrahedra as done in \cite{Suarez2021}, we also find finite orbits clustering around a different limiting cycle not yet formally identified. Interestingly, we find many orbits longer than the theoretical bound of $37$ introduced by Alder \cite{adler1983bisection} and analyzed with a different normalization procedure in \cite{PADRON2025555}. Overall, we have shown in Theorem \ref{thm:length} that there exist finite orbits of lengths $n\in\{3,4,7,8,13,14\} \cup \{21,\dots,26\}\cup\{36,37,38,39,42,43,44\}$ under the LEB.

However, there are also a number of tetrahedra whose orbits do not seem to be finite. We mention in this respect the cases of  the regular tetrahedron, cube corner tetrahedron, wedge tetrahedron, and needle tetrahedron, studied in the present work. Nevertheless, there is a good indication that their orbits may be still attracted to a certain closed region, yet to be formally identified. Our results provide a complementary approach to that presented in \cite{Suarez2021,Trujillo2024,Padron2023,PADRON2025555}. We use a different normalization procedure, but our results are highly compatible with theirs.

In the quest of proving the regularity of the tetrahedral partitions produced by the LEB refinements, we will develop our normalization procedure to provide a path for the developement of the  corresponding regularity proof in the spirit of the approach by Perdomo and Plaza \cite{Perdomo2013,Perdomo2014}. For such a proof to work, we need to endow the space of tetrahedron $\mathcal T$ with an appropriate metric. In Section \ref{sec:norm}, we  introduced a possible choice of such a metric, see Definition \ref{def:metric}. In the future work, we will explore the possibility to use that metric to identify closed regions in $\mathcal T$ and to eventually construct the desired regularity  proof for the LEB algorithm.

\begin{table}[ht]
\centering
\footnotesize
\begin{tabularx}{0.6\textwidth}{@{}lX@{}}
\hline
ID & Tetrahedron \\
\hline
$\tau_{1}$ &  \begin{minipage}{\linewidth}$\left(\frac{21}{76} + \frac{\sqrt{1079}}{76} i, \left(\frac{13}{19} + \frac{233 \sqrt{1079}}{41002} i, \frac{\sqrt{260424203}}{41002}\right)\right)$\end{minipage} \\
$\tau_{2}$ &  \begin{minipage}{\linewidth}$\left(\frac{53}{108} + \frac{\sqrt{1079}}{108} i, \left(\frac{61}{108} + \frac{7 \sqrt{1079}}{116532} i, \frac{\sqrt{82239222}}{19422}\right)\right)$\end{minipage} \\
$\tau_{3}$ &  \begin{minipage}{\linewidth}$\left(\frac{43}{88} + \frac{\sqrt{1759}}{88} i, \left(\frac{5}{11} + \frac{225 \sqrt{1759}}{19349} i, \frac{\sqrt{491580694}}{38698}\right)\right)$\end{minipage} \\
$\tau_{4}$ &  \begin{minipage}{\linewidth}$\left(\frac{20}{41} + \frac{3 \sqrt{33}}{41} i, \left(\frac{39}{82} + \frac{20 \sqrt{33}}{4059} i, \frac{\sqrt{17187159}}{8118}\right)\right)$\end{minipage} \\
$\tau_{5}$ &  \begin{minipage}{\linewidth}$\left(\frac{55}{116} + \frac{\sqrt{2543}}{116} i, \left(\frac{27}{58} + \frac{1299 \sqrt{2543}}{147494} i, \frac{\sqrt{936808141}}{73747}\right)\right)$\end{minipage} \\
$\tau_{6}$ &  \begin{minipage}{\linewidth}$\left(\frac{39}{80} + \frac{\sqrt{1759}}{80} i, \left(\frac{1}{2} + \frac{41 \sqrt{1759}}{3518} i, \frac{\sqrt{111722885}}{17590}\right)\right)$\end{minipage} \\
$\tau_{7}$ &  \begin{minipage}{\linewidth}$\left(\frac{5}{17} + \frac{2 \sqrt{15}}{17} i, \left(\frac{12}{17} + \frac{217 \sqrt{15}}{4080} i, \frac{\sqrt{3239265}}{4080}\right)\right)$\end{minipage} \\
$\tau_{8}$ &  \begin{minipage}{\linewidth}$\left(\frac{12}{25} + \frac{\sqrt{274}}{50} i, \left(\frac{59}{100}  - \frac{2 \sqrt{274}}{3425} i, \frac{\sqrt{1740311}}{2740}\right)\right)$\end{minipage} \\
$\tau_{9}$ &  \begin{minipage}{\linewidth}$\left(\frac{1}{2} + \frac{\sqrt{15}}{12} i, \left(\frac{41}{96} + \frac{\sqrt{15}}{240} i, \frac{\sqrt{63515}}{480}\right)\right)$\end{minipage} \\
$\tau_{10}$ &  \begin{minipage}{\linewidth}$\left(\frac{5}{17} + \frac{2 \sqrt{15}}{17} i, \left(\frac{10}{17} + \frac{263 \sqrt{15}}{4080} i, \frac{\sqrt{3239265}}{4080}\right)\right)$\end{minipage} \\
$\tau_{11}$ &  \begin{minipage}{\linewidth}$\left(\frac{39}{80} + \frac{\sqrt{1759}}{80} i, \left(\frac{19}{40} + \frac{939 \sqrt{1759}}{70360} i, \frac{\sqrt{111722885}}{17590}\right)\right)$\end{minipage} \\
$\tau_{12}$ &  \begin{minipage}{\linewidth}$\left(\frac{47}{96} + \frac{\sqrt{1055}}{96} i, \left(\frac{13}{32} + \frac{29 \sqrt{1055}}{33760} i, \frac{\sqrt{40204995}}{12660}\right)\right)$\end{minipage} \\
$\tau_{13}$ &  \begin{minipage}{\linewidth}$\left(\frac{39}{80} + \frac{\sqrt{2319}}{80} i, \left(\frac{1}{2} + \frac{55 \sqrt{2319}}{4638} i, \frac{\sqrt{147291285}}{23190}\right)\right)$\end{minipage} \\
$\tau_{14}$ &  \begin{minipage}{\linewidth}$\left(\frac{55}{116} + \frac{\sqrt{2543}}{116} i, \left(\frac{14}{29} + \frac{622 \sqrt{2543}}{73747} i, \frac{\sqrt{936808141}}{73747}\right)\right)$\end{minipage} \\
$\tau_{15}$ &  \begin{minipage}{\linewidth}$\left(\frac{57}{116} + \frac{\sqrt{2319}}{116} i, \left(\frac{27}{58} + \frac{1187 \sqrt{2319}}{134502} i, \frac{\sqrt{854289453}}{67251}\right)\right)$\end{minipage} \\
$\tau_{16}$ &  \begin{minipage}{\linewidth}$\left(\frac{1}{2} + \frac{\sqrt{15}}{12} i, \left(\frac{55}{96} + \frac{\sqrt{15}}{240} i, \frac{\sqrt{63515}}{480}\right)\right)$\end{minipage} \\
$\tau_{17}$ &  \begin{minipage}{\linewidth}$\left(\frac{5}{11} + \frac{\sqrt{274}}{22} i, \left(\frac{21}{44} + \frac{49 \sqrt{274}}{6028} i, \frac{\sqrt{19143421}}{6028}\right)\right)$\end{minipage} \\
$\tau_{18}$ &  \begin{minipage}{\linewidth}$\left(\frac{1}{2} + \frac{\sqrt{15}}{12} i, \left(\frac{41}{96}  - \frac{\sqrt{15}}{240} i, \frac{\sqrt{63515}}{480}\right)\right)$\end{minipage} \\
$\tau_{19}$ &  \begin{minipage}{\linewidth}$\left(\frac{23}{72} + \frac{\sqrt{1055}}{72} i, \left(\frac{7}{12} + \frac{97 \sqrt{1055}}{12660} i, \frac{\sqrt{26803330}}{12660}\right)\right)$\end{minipage} \\
$\tau_{20}$ &  \begin{minipage}{\linewidth}$\left(\frac{41}{86} + \frac{\sqrt{1243}}{86} i, \left(\frac{45}{86}  - \frac{39 \sqrt{1243}}{106898} i, \frac{\sqrt{678962647}}{53449}\right)\right)$\end{minipage} \\
$\tau_{21}$ &  \begin{minipage}{\linewidth}$\left(\frac{5}{17} + \frac{2 \sqrt{15}}{17} i, \left(\frac{41}{68} + \frac{287 \sqrt{15}}{4080} i, \frac{\sqrt{3239265}}{4080}\right)\right)$\end{minipage} \\
$\tau_{22}$ &  \begin{minipage}{\linewidth}$\left(\frac{6}{19} + \frac{\sqrt{274}}{38} i, \left(\frac{45}{76} + \frac{315 \sqrt{274}}{20824} i, \frac{\sqrt{66131818}}{20824}\right)\right)$\end{minipage} \\
$\tau_{23}$ &  \begin{minipage}{\linewidth}$\left(\frac{47}{96} + \frac{\sqrt{1055}}{96} i, \left(\frac{19}{32}  - \frac{29 \sqrt{1055}}{33760} i, \frac{\sqrt{40204995}}{12660}\right)\right)$\end{minipage} \\
$\tau_{24}$ &  \begin{minipage}{\linewidth}$\left(\frac{5}{17} + \frac{\sqrt{274}}{34} i, \left(\frac{10}{17} + \frac{331 \sqrt{274}}{18632} i, \frac{\sqrt{59170574}}{18632}\right)\right)$\end{minipage} \\
$\tau_{25}$ &  \begin{minipage}{\linewidth}$\left(\frac{39}{80} + \frac{\sqrt{2319}}{80} i, \left(\frac{19}{40} + \frac{1219 \sqrt{2319}}{92760} i, \frac{\sqrt{147291285}}{23190}\right)\right)$\end{minipage} \\
$\tau_{26}$ &  \begin{minipage}{\linewidth}$\left(\frac{41}{86} + \frac{\sqrt{1243}}{86} i, \left(\frac{41}{86} + \frac{39 \sqrt{1243}}{106898} i, \frac{\sqrt{678962647}}{53449}\right)\right)$\end{minipage} \\
$\tau_{27}$ &  \begin{minipage}{\linewidth}$\left(\frac{12}{25} + \frac{\sqrt{274}}{50} i, \left(\frac{41}{100} + \frac{2 \sqrt{274}}{3425} i, \frac{\sqrt{1740311}}{2740}\right)\right)$\end{minipage} \\
$\tau_{28}$ &  \begin{minipage}{\linewidth}$\left(\frac{53}{108} + \frac{\sqrt{1079}}{108} i, \left(\frac{47}{108}  - \frac{7 \sqrt{1079}}{116532} i, \frac{\sqrt{82239222}}{19422}\right)\right)$\end{minipage} \\
$\tau_{29}$ &  \begin{minipage}{\linewidth}$\left(\frac{5}{17} + \frac{2 \sqrt{15}}{17} i, \left(\frac{47}{68} + \frac{193 \sqrt{15}}{4080} i, \frac{\sqrt{3239265}}{4080}\right)\right)$\end{minipage} \\
$\tau_{30}$ &  \begin{minipage}{\linewidth}$\left(\frac{41}{88} + \frac{\sqrt{2543}}{88} i, \left(\frac{21}{44} + \frac{1251 \sqrt{2543}}{111892} i, \frac{\sqrt{710682038}}{55946}\right)\right)$\end{minipage} \\
$\tau_{31}$ &  \begin{minipage}{\linewidth}$\left(\frac{57}{116} + \frac{\sqrt{2319}}{116} i, \left(\frac{15}{29} + \frac{566 \sqrt{2319}}{67251} i, \frac{\sqrt{854289453}}{67251}\right)\right)$\end{minipage} \\
$\tau_{32}$ &  \begin{minipage}{\linewidth}$\left(\frac{43}{88} + \frac{\sqrt{1759}}{88} i, \left(\frac{23}{44} + \frac{859 \sqrt{1759}}{77396} i, \frac{\sqrt{491580694}}{38698}\right)\right)$\end{minipage} \\
$\tau_{33}$ &  \begin{minipage}{\linewidth}$\left(\frac{20}{41} + \frac{3 \sqrt{33}}{41} i, \left(\frac{43}{82}  - \frac{20 \sqrt{33}}{4059} i, \frac{\sqrt{17187159}}{8118}\right)\right)$\end{minipage} \\
$\tau_{34}$ &  \begin{minipage}{\linewidth}$\left(\frac{19}{72} + \frac{\sqrt{1079}}{72} i, \left(\frac{49}{72} + \frac{473 \sqrt{1079}}{77688} i, \frac{\sqrt{27413074}}{12948}\right)\right)$\end{minipage} \\
$\tau_{35}$ &  \begin{minipage}{\linewidth}$\left(\frac{41}{88} + \frac{\sqrt{2543}}{88} i, \left(\frac{5}{11} + \frac{323 \sqrt{2543}}{27973} i, \frac{\sqrt{710682038}}{55946}\right)\right)$\end{minipage} \\
$\tau_{36}$ &  \begin{minipage}{\linewidth}$\left(\frac{1}{2} + \frac{\sqrt{15}}{12} i, \left(\frac{55}{96}  - \frac{\sqrt{15}}{240} i, \frac{\sqrt{63515}}{480}\right)\right)$\end{minipage} \\
$\tau_{37}$ &  \begin{minipage}{\linewidth}$\left(\frac{19}{72} + \frac{\sqrt{1079}}{72} i, \left(\frac{7}{12} + \frac{101 \sqrt{1079}}{12948} i, \frac{\sqrt{27413074}}{12948}\right)\right)$\end{minipage} \\
$\tau_{38}$ &  \begin{minipage}{\linewidth}$\left(\frac{21}{76} + \frac{\sqrt{1079}}{76} i, \left(\frac{45}{76} + \frac{613 \sqrt{1079}}{82004} i, \frac{\sqrt{260424203}}{41002}\right)\right)$\end{minipage} \\
$\tau_{39}$ &  \begin{minipage}{\linewidth}$\left(\frac{21}{68} + \frac{\sqrt{1055}}{68} i, \left(\frac{41}{68} + \frac{669 \sqrt{1055}}{71740} i, \frac{\sqrt{227828305}}{35870}\right)\right)$\end{minipage} \\
\hline
\end{tabularx}
\caption{List of tetrahedra in the orbit of T1 defined in Table \ref{tab:T1-6}.}
\label{tab:tet_nodes_39}
\end{table}

\bibliographystyle{elsarticle-num}
\bibliography{references}
\end{document}